\documentclass[twoside]{article}
\usepackage{amsmath,amsthm,amssymb,stmaryrd,latexsym,amsbsy,xy}
\usepackage{amscd, amssymb, verbatim}
\usepackage{pstricks,pst-plot,pst-node}
\usepackage{amsfonts}
\usepackage{anysize}
\usepackage[pdftex,bookmarks=true]{hyperref}
\usepackage{fancyhdr}

\marginsize{3.3cm}{3.3cm}{3.315cm}{3.315cm}

\newtheorem{theorem}{Theorem}[section]

\newtheorem{lemma}[theorem]{Lemma}
\newtheorem{prop}[theorem]{Proposition}

\newtheorem{conj}[theorem]{Conjecture}
\newtheorem{conj/cor}[theorem]{``Corollary''}
\newtheorem{conj/prop}[theorem]{``Proposition''}
\newtheorem{conj/theorem}[theorem]{``Theorem''}

\numberwithin{equation}{section}
\newcommand{\bpf}{\noindent {\em Proof.  }}
\newcommand{\bsp}{\noindent {\em Sketch of proof.  }}
\newcommand{\epf}{\qed \vspace{+10pt}}

\newcommand{\ep}{\qed \vspace{+10pt}}

\newcommand{\Hilb}{\mathrm{Hilb}}
\newcommand{\Sym}{\mathrm{Sym}}
\newcommand{\Ast}{A^*}
\newcommand{\AT}{A_{\bT}^*}

\newcommand{\Aorb}{A_{\mathrm{orb}}^*}
\newcommand{\ATorb}{A_{\bT, \mathrm{orb}}^*}
\newcommand{\HTorb}{H_{\bT, \mathrm{orb}}^*}
\newcommand{\QAT}{QA_{\bT}^*}
\newcommand{\QATorb}{QA_{\bT, \mathrm{orb}}^*}

\newcommand{\orb}{\mathrm{orb}}
\newcommand{\vir}{\mathrm{vir}}
\newcommand{\mov}{\mathrm{mov}}

\newcommand{\HC}{\mathrm{HC}}
\newcommand{\Spec}{\mathrm{Spec}}

\newcommand{\lcm}{\mathrm{lcm}}

\newcommand{\M}{\overline{M}}
\newcommand{\Mc}{\overline{M}^{\circ}}

\newcommand{\llangle}{\langle \langle}
\newcommand{\rrangle}{\rangle \rangle}

\newcommand{\I}{\overline{I}}
\newcommand{\age}{\mathrm{age}}
\newcommand{\Aut}{\mathrm{Aut}}
\newcommand{\conn}{\mathrm{conn}}

\newcommand{\stab}{\mathrm{Stab}}
\newcommand{\eT}{e_{\mathbb{T}}}
\newcommand{\ev}{\mathrm{ev}}
\newcommand{\MT}{\M_{\bT}}
\newcommand{\spa}{\textrm{ }}
\newcommand{\Aone}{\mathcal{A}_1}
\newcommand{\Ar}{\mathcal{A}_{r}}
\newcommand{\An}{\mathcal{A}_{n}}
\newcommand{\tC}{\tilde{C}}
\newcommand{\tf}{\tilde{f}}
\newcommand{\tD}{\tilde{D}}
\newcommand{\Eij}{\cE_{ij}}
\newcommand{\ttt}{(t_1+t_2)}
\newcommand{\qorb}{*_{\orb}}
\newcommand{\qcrp}{*_{\mathrm{crep}}}
\newcommand{\GW}{\mathrm{GW}}
\newcommand{\Zp}{\mathrm{Z}^\prime}

\newcommand{\fa}{\mathfrak{a}}
\newcommand{\fB}{\mathfrak{B}}

\newcommand{\fm}{\mathfrak{m}}
\newcommand{\fp}{\mathfrak{p}}
\newcommand{\fS}{\mathfrak{S}}

\newcommand{\cB}{\mathcal{B}}
\newcommand{\cC}{\mathcal{C}}
\newcommand{\cD}{\mathcal{D}}
\newcommand{\cE}{\mathcal{E}}
\newcommand{\cF}{\mathcal{F}}
\newcommand{\cH}{\mathcal{H}}
\newcommand{\cI}{\mathcal{I}}
\newcommand{\cL}{\mathcal{L}}

\newcommand{\cO}{\mathcal{O}}
\newcommand{\cP}{\mathcal{P}}
\newcommand{\cQ}{\mathcal{Q}}

\newcommand{\cV}{\mathcal{V}}
\newcommand{\cX}{\mathcal{X}}

\newcommand{\bC}{\mathbb{C}}
\newcommand{\bH}{\mathbb{H}}
\newcommand{\bN}{\mathbb{N}}
\newcommand{\bP}{\mathbb{P}}
\newcommand{\bQ}{\mathbb{Q}}
\newcommand{\bT}{\mathbb{T}}
\newcommand{\bZ}{\mathbb{Z}}

\newcommand{\al}{\alpha}
\newcommand{\be}{\beta}
\newcommand{\g}{\gamma}

\newcommand{\de}{\delta}
\newcommand{\si}{\sigma}
\newcommand{\Si}{\Sigma}
\newcommand{\la}{\lambda}
\newcommand{\La}{\Lambda}

\newcommand{\wt}{\widetilde}

\begin{document}

\title{\Large \bf Orbifold quantum cohomology of the symmetric product of $\Ar$}

\author{Wan Keng Cheong\footnote{ Department of Mathematics, California Institute of Technology, Pasadena, CA, USA. \newline Email address: keng@caltech.edu}}
\date{}
\maketitle

\pagenumbering{roman}
\thispagestyle{empty}
\begin{abstract}
Let $\Ar$ be the minimal resolution of the type $A_r$ surface singularity.
We study the equivariant quantum cohomology ring of the $n$-fold symmetric product stack $[\Sym^n(\Ar)]$ of $\Ar$. We calculate the operators of quantum multiplication by divisor classes. Under the assumption of the nonderogatory conjecture, these operators completely determine the ring structure, which gives an affirmative answer to the Crepant Resolution Conjecture on $[\Sym^n(\Ar)]$ and $\Hilb^n(\Ar)$. More strikingly, this allows us to complete a tetrahedron of equivalences relating the Gromov-Witten theories of $[\Sym^n(\Ar)]$/$\Hilb^n(\Ar)$ and the relative Gromov-Witten/Donaldson-Thomas theories of $\Ar \times \bP^1$.
\end{abstract}

\tableofcontents

\newpage 

\pagenumbering{arabic}

\setcounter{section}{-1}
\section{Overview}
\subsection{Results}
In physics, it is believed that the string theory on a quotient space and the string theory on any crepant resolution should belong to the same family. Over the years, this principle has been put into various mathematical frameworks. Among them, we are particularly interested in the formulations pioneered by Ruan in the context of Gromov-Witten theory (see e.g. \cite{R, BG, CoIT, CoR}). 

Let $\Ar$ be the minimal resolution of the type $A_r$ surface singularity. The symmetric group $\fS_n$ acts on the $n$-fold Cartesian product $\Ar^n$ by permuting coordinates.
Thus, we obtain a quotient scheme $\Sym^n(\Ar):=\Ar^n/\fS_n$, the $n$-fold symmetric product of $\Ar$, and a quotient stack $[\Sym^n(\Ar)]$, the $n$-fold symmetric product stack of $\Ar$. The stack $[\Sym^n(\Ar)]$ is a smooth orbifold, whose coarse moduli space is none other than the symmetric product $\Sym^n(\Ar)$.

In this article, we compare the equivariant orbifold Gromov-Witten theory of the symmetric products of $\Ar$ with the equivariant Gromov-Witten theory of the crepant resolutions in the spirit of Bryan-Graber's Crepant Resolution Conjecture \cite{BG}. 

Let $T=\bC^{\times}\times \bC^{\times}$ be a two-dimensional torus. The (localized) $\bT$-equivariant cohomology of a point is generated by $t_1$ and $t_2$. Our main objects are the $3$-point functions
$$
\llangle \al_1, \al_2, \al_3 \rrangle^{[\Sym^n(\Ar)]}\in \bQ(t_1,t_2)[[u, s_1, \ldots, s_r]]
$$
which encode $3$-point extended Gromov-Witten invariants of $[\Sym^n(\Ar)]$ (see \eqref{k-point ext}). These generating functions add a multiplicative structure to the equivariant Chen-Ruan cohomology $\HTorb([\Sym^n(\Ar)]; \bQ)$. The multiplication so obtained is called the small orbifold quantum product. 

The quotient space $\Sym^n(\Ar)$ admits a unique crepant resolution of singularities, namely the Hilbert scheme $\Hilb^n(\Ar)$ of $n$ points in $\Ar$. The $\bT$-equivariant quantum cohomology of $\Hilb^n(\Ar)$ has been explored by Maulik and Oblomkov in \cite{MO1}, so we need only deal with the quantum ring of the orbifold $[\Sym^n(\Ar)]$.
We fully cover $2$-point extended Gromov-Witten invariants of $[\Sym^n(\Ar)]$ and find that the calculation of these invariants is tantamount to the question of counting certain branched covers of rational curves. Our discovery can be summarized in the following statement.
\begin{theorem}\label{two-point theorem}
Two-point extended equivariant Gromov-Witten invariants of $[\Sym^n(\Ar)]$ are expressible in terms of equivariant orbifold Poincar\'e pairings and one-part double Hurwitz numbers. 
\end{theorem}
One-part double Hurwitz numbers, as shown by Goulden, Jackson and Vakil \cite{GJV}, admit explicit closed formulas (c.f. \eqref{GJV}) and therefore Theorem \ref{two-point theorem} provides a complete solution to the divisor operators, i.e. the operators of quantum multiplication by divisor classes. These operators correspond naturally to the divisor operators on the Hilbert scheme $\Hilb^n(\Ar)$:
\begin{theorem}\label{hilbert}
After making the change of variables $q=-e^{iu}$, where $i^2=-1$, and extending scalars to an appropriate field $F$, there is a linear isomorphism of equivariant quantum cohomologies
$$L: \HTorb([\Sym^n(\Ar)]; F) \to H_{\bT}^*(\Hilb^n(\Ar); F)$$
which preserves gradings, Poincar\'e pairings and respects small quantum product by divisors. In other words, for any Chen-Ruan cohomology classes $\al_1, \al_2$ and divisor $D$, we have the following identity for $3$-point functions:
$$
\llangle \al_1, D, \al_2 \rrangle^{[\Sym^n(\Ar)]}=\langle L(\al_1), L(D), L(\al_2) \rangle^{\Hilb^n(\Ar)}.
$$
Here $\langle -,-,- \rangle^{\Hilb^n(\Ar)}$ are the $3$-point functions of $\Hilb^n(\Ar)$ in variables $t_1, t_2$, $q, s_1, \ldots, s_r$ (see \eqref{k-point hilb}).
\end{theorem}

In addition to the relation to the Hilbert schemes, the orbifold theory is in connection with the relative Gromov-Witten theory of threefolds:
\begin{theorem}\label{relative}
Given cohomology-weighted partitions $\la_1(\vec{\eta}_1)$, $\la_2(\vec{\eta}_2)$ of $n$ and $\al=1(1)^n$, $(2)$ or $D_k$, $k=1,\ldots,r$  (see Section \ref{description} and Section \ref{relative GW} for the corresponding classes), we have
$$
\llangle \la_1(\vec{\eta}_1), \al, \la_2(\vec{\eta}_2)\rrangle^{[\Sym^n(\Ar)]}=\GW (\Ar \times \bP^1)_{\la_1(\vec{\eta}_1), \al, \la_2(\vec{\eta}_2)},
$$
where the right hand side is a shifted partition function (c.f. \eqref{shifted partition}).
\end{theorem}

\subsection{Equivalence with other theories}\label{tetrahedron}
The above theorems form a triangle of equivalences. We can include the Donaldson-Thomas theory to make up a tetrahedron. In fact, Theorem \ref{hilbert} and \ref{relative}, in conjunction with the results of \cite{M, MO1, MO2}, establish the following equivalences for divisor operators.

\vspace{2.8cm}
\begin{center}
\begin{picture}(80,50)(-30,-55)
\put(-50,0){\line(1,1){50}}
\multiput(-50,0)(8,0){13}{\line(1,0){2}}
\put(50,0){\line(-1,1){50}}

\put(-50,0){\line(4,-1){72.5}}
\put(50,0){\line(-3,-2){27}}

\put(0,50){\line(1,-3){22.7}}

\put(0,67){\makebox(0,0){Orbifold quantum }}
\put(0,57){\makebox(0,0){cohomology of $[\Sym(\Ar)]$ }}

\put(-95,10){\makebox(0,0){Relative}}
\put(-95,0){\makebox(0,0){Gromov-Witten}}
\put(-95,-10){\makebox(0,0){theory of $\Ar \times \bP^1$}}

\put(30,-25){\makebox(0,0){Quantum cohomology }}
\put(30,-35){\makebox(0,0){of $\Hilb(\Ar)$ }}

\put(100,10){\makebox(0,0){Relative}}
\put(100,0){\makebox(0,0){Donaldson-Thomas}}
\put(100,-10){\makebox(0,0){theory of $\Ar \times \bP^1$}}
\put(0,-50){\makebox(0,0){\small{Figure 1: A tetrahedron of equivalences. }}}
\end{picture}
\end{center}

The nonderogatory conjecture (Conjecture \ref{nonderogatory}) of Maulik and Oblomkov implies that the quantum ring will be generated by divisor classes, and so the linear isomorphism $L$ in Theorem \ref{hilbert} will be a ring isomorphism. 

Before the study of the Gromov-Witten theory of $[\Sym^n(\Ar)]$, the case of the affine plane $\bC^2$ was the only known example for the above tetrahedron to hold for all operators (c.f. \cite{BG, BP, OP1, OP2}). If the nonderogatory conjecture is assumed, these four theories will be equivalent in our case of $\Ar$ as well. The base triangle of ``equivalences'' is the work of Maulik and Oblomkov. And the triangle facing the rightmost corner is worked out in this paper:
\begin{conj/prop}\label{triangle}
Let $L$ be as in Theorem \ref{hilbert} and $\la_1(\vec{\eta}_1)$, $\la_2(\vec{\eta}_2)$, $\la_3(\vec{\eta}_3)$ any cohomology-weighted partitions of $n$. Assuming the nonderogatory conjecture, the identities
\begin{eqnarray*}
\llangle \la_1(\vec{\eta}_1), \la_2(\vec{\eta}_2), \la_3(\vec{\eta}_3) \rrangle^{[\Sym^n(\Ar)]} &=& \langle L(\la_1(\vec{\eta}_1)), L(\la_2(\vec{\eta}_2)), L(\la_3(\vec{\eta}_3)) \rangle^{\Hilb^n(\Ar)}\\
&=& \GW(\Ar \times \bP^1)_{\la_1(\vec{\eta}_1), \la_2(\vec{\eta}_2), \la_3(\vec{\eta}_3)}
\end{eqnarray*}
hold under the substitution $q=-e^{iu}$.
\end{conj/prop}

Once the WDVV equations are used, we can make a more general statement on $[\Sym^n(\Ar)]$ and $\Hilb^n(\Ar)$.
\begin{conj/prop}\label{CRC}
Let $q=-e^{iu}$. Assuming the nonderogatory conjecture, the map $L$, in Theorem \ref{hilbert}, equates the extended multipoint functions of $[\Sym^n(\Ar)]$ to the multipoint functions of $\Hilb^n(\Ar)$. Moreover, these functions are rational functions in $t_1, t_2$, $q, s_1, \ldots, s_r$.
(Multipoint functions are those with at least three insertions). 
\end{conj/prop}
This answers positively the Crepant Resolution Conjecture, proposed by Bryan and Graber, on the symmetric product case. We will see that ``Proposition'' \ref{triangle} and \ref{CRC} are valid in the case of $n=2, r=1$ even without presuming the nonderogatory conjecture (c.f. Section \ref{s:example}). 

\subsection{Outline of the paper} 
The aim of Section \ref{s:1} is to give a brief introduction to the resolved surface $\Ar$ and collect some basic facts required for the paper.

In Section \ref{s:2}, we recall Chen-Ruan's orbifold cohomology for a symmetric product and construct certain bases for the $\bT$-equivariant orbifold cohomology. We also present an algorithm to express these bases in terms of $\bT$-fixed point basis.

Section \ref{s:3} is devoted to reviewing some background on orbifold Gromov-Witten theory and defining extended Gromov-Witten invariants and their connected counterparts.

Section \ref{s:4} is the main theme. We calculate $2$-point extended invariants of nonzero degrees by virtual localization. The results we obtain prove Theorem \ref{two-point theorem} and allow a combinatorial description of any divisor operator for $[\Sym^n(\Ar)]$.

In Section \ref{s:5}, we prove Theorem \ref{hilbert} and \ref{relative}. We also exhibit a simple example ($n=2, r=1$) for which the correspondence in Theorem \ref{hilbert} is actually an isomorphism of quantum rings. 

In Section \ref{s:6}, a certain nondegeneracy conjecture is described. Under the assumption of this conjecture, an equivalence between the equivariant Gromov-Witten theories of the symmetric product stack and the Hilbert scheme of points will be obtained and so will ``Proposition'' \ref{triangle}. In the end, we discuss multipoint functions of $[\Sym^n(\Ar)]$ and the full version of the Crepant Resolution Conjecture (``Proposition'' \ref{CRC}).

\subsection{Notation and convention}\label{notation}
The following notations will be used without further comment. Some other notations will be introduced along the way.

\begin{enumerate}
\item To avoid doubling indices, we identify 
$$A^i(X)=H^{2i}(X; \bQ), \spa A_i(X)=H_{2i}(X; \bQ) \spa \textrm{and } A_i(X; \bZ)=H_{2i}(X; \bZ),$$ 
just to name a few, for any complex variety $X$ to appear in this article (note that we drop $\bQ$ but not $\bZ$). They will be referred to as cohomology or homology groups rather than Chow groups. 

\item An orbifold $\cX$ is a smooth Deligne-Mumford stack of finite type over $\bC$. Denote by $c: \cX \to X$ the canonical map to the coarse moduli space.

\item For any positive integer $s$, $\mu_s$ is the cyclic subgroup of $\bC^\times$ of order $s$.
 
\item For any finite group $G$, $\cB G$ is the classifying stack of $G$, i.e. $[\Spec \spa \bC/G]$.

\item
 \begin{enumerate}
 \item $\bT=(\bC^{\times})^2$ is always a two-dimensional torus.

 \item $t_1$, $t_2$ are the generators of the $\bT$-equivariant cohomology $A_{\bT}^*(\mathrm{point})$ of a point, that is, $A_{\bT}^*(\mathrm{point})=\bQ[t_1,t_2].$

 \item $V_{\fm}=V\otimes_{\bQ[t_1,t_2]}\bQ(t_1,t_2)$ for each $\bQ[t_1,t_2]$-module $V$.
 \end{enumerate}

\item Given any object $\cO$, $\cO^n$ means that $\cO$ repeats itself $n$ times.

\item For $i=1, 2$, $\epsilon_i$ is a function on the set of nonnegative integers such that
\begin{equation*}
\ \epsilon_i(m)=
\begin{cases}
0 &\text{if $m < i$;}\\
1 &\text{if $m \geq i$.}
\end{cases}
\end{equation*}

\item Given a partition $\si$ of a nonnegative integer. 
 \begin{enumerate}
 \item $\ell(\si)$ is the length of $\si$. 
 
 \item Unless otherwise stated, $\si$ is presumed to be written as
 $$\si=(\si_1, \ldots, \si_{\ell(\si)}).$$ 
  To make a emphasis, if $\si_k$ is another partition, it is simply $(\si_{k1}, \ldots, \si_{k\ell(\si_k)})$.
  
 \item $|\si|=n$ if $\si_1+\cdots +\si_{\ell(\si)}=n$.

 \item Let $\vec{\al}:=(\al_1, \ldots, \al_{\ell(\si)})$ be an $\ell(\si)$-tuple of cohomology classes associated to $\si$ so that we may form a cohomology-weighted partition $\si(\vec{\al}):=\si_1(\al_1)\cdots \si_{\ell(\si)}(\al_{\ell(\si)})$. The group $\Aut(\si(\vec{\al}))$ is defined to be the group of permutations on $\{1, 2, \ldots, \ell(\si) \}$ fixing
$$( \spa (\si_1, \al_1),\ldots, (\si_{\ell(\si)}, \al_{\ell(\si)})\spa ).$$
Let $\Aut(\si)$ be the group $\Aut(\si(\vec{\al}))$ when all entries of $\vec{\al}$ are identical. Its order is simply $\prod_{i=1}^n m_i!$ if $\si=(1^{m_1}, \ldots, n^{m_n})$.
 \item $o(\si)=\lcm(|\si_1|, \ldots, |\si_{\ell(\si)}|)$ is the order of any permutation of cycle type $\si$.
 \item $(2):=(1^{n-2}, 2)$ and $1:=(1^n)$ are partitions of length $n-1$ and length $n$ respectively.
 \end{enumerate}
 
\end{enumerate}

\section{Resolutions of cyclic quotient surface singularities} \label{s:1}
We fix a positive integer $r$ once and for all. Let the cyclic group $\mu_{r+1}$ act on $\bC^2$ by the diagonal matrices
$$
\left(\begin{matrix}
\zeta & 0\\
0 & \zeta^{-1}
\end{matrix}\right),\\
$$ 
where $\zeta \in \mu_{r+1}$. The quotient $\bC^2/\mu_{r+1}$ is a surface singularity. We denote by
$$\pi: \Ar \to \bC^2/\mu_{r+1}$$
its minimal resolution. It is actually well-known that $\pi$ can be obtained via a sequence of $\lfloor \frac{r+1}{2} \rfloor$ blow-ups at the unique singularity.
The exceptional locus $\mathrm{Ex}(\pi)$ of $\pi$ is a chain of $(-2)$-curves,
$\bigcup_{i=1}^r E_i,$
with $E_{i-1}$ and $E_i$ intersect transversally. 
The intersection numbers of the exceptional curves are given by
\[
E_i \cdot E_j =
\begin{cases}
-2 &\text{if $i=j$;}\\
1 &\text{if $|i-j|=1$;}\\
0 &\text{otherwise.}
\end{cases}
\]
In particular, the intersection matrix is negative definite (as expected from the general theory of complex surfaces). Additionally, $E_1, \ldots, E_r$ give a basis for $A_1(\Ar; \bZ)$. 
We also have two noncompact curves $E_0$ and $E_{r+1}$ attached to $E_1$ and $E_r$ respectively. $E_0$ (resp. $E_{r+1}$) can be arranged to map to the $\mu_{r+1}$-orbit of $x$-axis (resp. $y$-axis).

The natural action of $\bT$ on $\bC^2$ comes with tangent weights $t_1$ and $t_2$ at the origin. It commutes with the $\mu_{r+1}$-action, so we have an induced $\bT$-action on the quotient $\bC^2/\mu_{r+1}$ and thus on the resolved surface $\Ar$. We fix these actions of $\bT$ throughout the article.

The $\bT$-invariant curves on $\Ar$ are $E_1, \ldots, E_r$. The $\bT$-fixed points are the points at the nodes of the chain $\cup_{i=0}^{r+1}E_i$ of curves. Precisely, they are
$$x_1, \ldots, x_{r+1},$$
where $\{ x_i \}=E_{i-1} \cap E_i$. Let's assume that $L_i$ and $R_i$ are respectively the weights of the $\bT$-action on the tangent spaces to $E_{i-1}$ and $E_i$ at $x_i$. We have $L_1=(r+1)t_1$, $R_{r+1}=(r+1)t_2$ and the following equalities
$$
L_i+R_i=t_1+t_2, \spa R_i=-L_{i+1},
$$
for each $i=1, \ldots r$.

\begin{center}
\begin{picture}(50,60)(40,-20)
\put(-105,20){$E_0$}

\put(-56,0){\line(-4,3){50}}
\put(-56,0){\line(4,3){35}}
\put(-56,0){\vector(-4,3){12}}
\put(-126,0){$(r+1)t_1=L_1$}
\put(-56,0){\vector(4,3){12}}
\put(-47,0){$R_1$}
\put(-58,-2){$\bullet$}
\put(-58,-15){$x_1$}
\put(-56,0){\line(4,-3){10}}
\put(-56,0){\line(-4,-3){10}}

\put(-48,18){$L_2$}
\put(-28,21){\vector(-4,-3){12}}
\put(-30,19){$\bullet$}
\put(-32,27){$x_2$}

\put(0,0){\line(4,3){35}}
\put(0,0){\line(-4,3){35}}
\put(0,0){\line(-4,-3){10}}
\put(0,0){\line(4,-3){10}}

\put(42,10){$\cdot$}
\put(56,10){$\cdot$}
\put(70,10){$\cdot$}

\put(112,0){\line(-4,3){35}}
\put(112,0){\line(4,3){35}}
\put(112,0){\line(4,-3){10}}
\put(112,0){\line(-4,-3){10}}

\put(138,19){$\bullet$}
\put(136,27){$x_r$}
\put(140,21){\vector(4,-3){12}}
\put(147,18){$R_r$}

\put(168,0){\line(-4,3){35}}

\put(168,0){\vector(-4,3){12}}
\put(140,0){$L_{r+1}$}

\put(168,0){\vector(4,3){12}}
\put(178,0){$R_{r+1}=(r+1)t_2$}

\put(166,-2){$\bullet$}
\put(160,-15){$x_{r+1}$}
\put(168,0){\line(4,-3){10}}
\put(168,0){\line(-4,-3){10}}
\put(168,0){\line(4,3){50}}

\put(210,20){$E_{r+1}$}
\end{picture}
\end{center}
\begin{quote}\small{
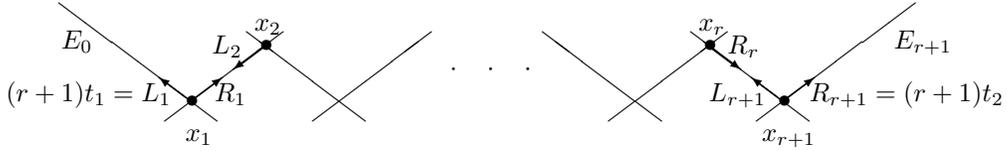
Figure 2: The middle chain is the exceptional locus $\mathrm{Ex}(\pi)$. The labelled vectors stand for the tangent weights at the fixed points.}
\end{quote}

The above information will be sufficient for our calculation of Gromov-Witten invariants. Certainly, one can also compute explicitly to obtain
$$
(L_i, R_i)=((r-i+2)t_1+(1-i)t_2, \spa (-r+i-1)t_1+it_2).
$$

\section{Chen-Ruan cohomology} \label{s:2}
\subsection{Inertia stack}
Given any finite set $N$, let $\fS_N$ be the symmetric group on $N$ and 
$$X^N=\{(x_i)_{i\in N}: x_i\textrm{'s are elements of }X\}$$ 
is a set of $|N|$-tuples of elements of $X$.
We denote by $\fS_n$ the group $\fS_{\{1,\ldots, n\}}$ and by $X^n$ the set $X^{\{1,\ldots, n\}}$.

The symmetric group $\fS_n$ acts on $X^n$ by $g(z)_i=z_{g(i)}$, $\forall g\in \fS_n$, $z \in X^n.$ The $n$-fold symmetric product $\Sym^n(X)$ is defined to be $X^n/\fS_n$, and the $n$-fold symmetric product stack $[\Sym^n(X)]$ is defined to be the quotient stack $[X^n/\fS_n]$. The space $\Sym^n(X)$ is in general singular and is the coarse moduli space of the (smooth) orbifold $[\Sym^n(X)]$. 

There is a natural stack associated to the symmetric product, i.e. the inertia stack
$$I[\Sym^n(X)]:=\coprod_{s\in \bN}HomRep(\cB \mu_s, [\Sym^n(X)]),$$
where $HomRep(\cB \mu_s, [\Sym^n(X)])$ is the stack of representable morphisms from the classifying stack $\cB \mu_s$ to $[\Sym^n(X)]$. Moreover, $I[\Sym^n(X)]$ is isomorphic to the disjoint union of orbifolds
\begin{equation}\label{sectors}
\coprod_{[g]\in C}[X^n_g/C(g)],
\end{equation}
where $C$ is the set of conjugacy classes, $X^n_g$ is the $g$-fixed locus of $X^n$, and $C(g)$ is the centralizer of $g$.
The component $[X^n/\fS_n]$ is called the untwisted sector while all other components are called twisted sectors. As there is a one-to-one correspondence between the conjugacy classes of $\fS_n$ and the partitions of $n$, these sectors can be labelled with the partitions of $n$. If $[g]$ is the conjugacy class corresponds to the partition $\la$ and $\overline{C(g)}:=C(g)/\langle g \rangle$, we may write
$$ 
X(\la):=X^n_g/C(g), \textrm{ and }\overline{X(\la)}:=X^n_g/\overline{C(g)} \textrm{ (see below).}
$$

The Chen-Ruan cohomology
$$\Aorb([\Sym^n(X)])$$ 
is by definition the cohomology $\Ast(I[\Sym^n(X)])$ of the inertia stack (\cite{ChR1}).
By \eqref{sectors}, it is 
$\bigoplus_{[g]\in C}\Ast(X^n_g/C(g))=\bigoplus_{[g]\in C}\Ast(X^n_g)^{C(g)}.$ (For any orbifold $\mathcal{Y}$ with coarse moduli space $Y$, we identify $\Ast(\mathcal{Y})=\Ast(Y)$ by the pushforward $c_*:\Ast(\mathcal{Y}) \to \Ast(Y)$ defined by $c_*([\cV])=\frac{1}{s}[c(\cV)]$, where $\cV$ is a closed integral substack and $s$ is the order of the stabilizer of a generic geometric point of $\cV$).

The age (or the degree shifting number) of the sector $[X(\la)]$ is given by
$$\age(\la):=n-\ell(\la).$$
Additionally, the Chen-Ruan cohomology is graded by ages. If $\al\in A^i(X(\la))$, the orbifold (Chow) degree of $\al$ is defined to be $i+\age(\la)$. In other words,
$$A_{\mathrm{orb}}^*([\Sym^n(X)])=\bigoplus_{|\la|=n} A^{*-\age(\la)}(X(\la)).$$

When $X$ admits a $\bT$-action, we can see easily that there are induced $\bT$-actions on the spaces $X^n_g/C(g)$ ($\forall g\in \fS_n$) and $I[\Sym^n(X)]$. So we may put the above cohomologies into an equivariant context by considering $\bT$-equivariant cohomologies.

We may rigidify the inertia stack to remove the actions of $\mu_s$'s. Each $\cB \mu_s$ acts on the stack $HomRep(\cB \mu_s, [\Sym^n(X)])$ and the quotient by this action is a stack of gerbes banded by $\mu_s$ to $[\Sym^n(X)]$. The stack
$$
\I[\Sym^n(X)]:=\coprod_{s\in \bN}HomRep(\cB \mu_s, [\Sym^n(X)])/\cB \mu_s,
$$
is called the rigidified inertia stack of $[\Sym^n(X)]$. One of the reasons why we mention this is that the rigidified stack is where the evaluation maps land (see \eqref{ev}). 

For more details on the rigidification procedure, consult \cite{ACV, AGV1, AGV2}.
In fact, the procedure amounts to removing the action of the permutation $g$ from \eqref{sectors} -- $g$ acts trivially on $X^n_g$. Concretely, the stack $\I[\Sym^n(X)]$ is the disjoint union
$\coprod_{[g]\in C}[X^n_g/\overline{C(g)}]$ (or $\coprod_{|\la|=n}[\overline{X(\la)}]$). Its coarse moduli space $\coprod_{[g]\in C} X^n_g/C(g)$ is identical to that of the inertia stack. 

\subsection{Bases}
\subsubsection{A description}\label{description}
Let $X$ be a smooth toric surface. We need to understand the module structure of equivariant Chen-Ruan cohomology $\ATorb([\Sym^n(X)])$; particularly, we need a precise description of bases.
Given a partition $\la$ of $n$, we would like to give a basis for the cohomology $\AT(X^n_{g})^{C(g)}$, where $g\in \fS_n$ has cycle type $\la$.

The permutation $g$ has a cycle decomposition, i.e. a product of disjoint cycles (including 1-cycles),
$$g=g_1 \ldots g_{\ell(\la)}$$ 
with $g_i$ being a $\la_i$-cycle.
For each $i$, let $N_i$ be the minimal subset of $\{1,\ldots, n\}$ such that $g_i \in \fS_{N_i}$. Thus
$|N_i|=\la_i$ and $\coprod_{i=1}^{\ell(\la)}  N_i=\{1,\ldots, n\}.$ It is clear that
$$
X^n_g=\prod_{i=1}^{\ell(\la)} X^{N_i}_{g_{i}}, \textrm{ and } X^{N_i}_{g_i}\cong X.
$$

To the partition $\la$, we associate a $\ell(\la)$-tuple $\vec{\eta}=(\eta_1 \ldots \eta_{\ell(\la)})$ with entries in $\AT(X)$. Let's put
\begin{equation}\label{basis}
g(\vec{\eta}) =(|\Aut(\la(\vec{\eta}))|\prod_{i=1}^{\ell(\la)} \la_i)^{-1} \sum_{h\in C(g)}\bigotimes_{i=1}^{\ell(\la)} g_i^h(\eta_i) \in \AT(X^n_{g})^{C(g)}.
\end{equation}
This requires some explanations:
\begin{itemize}
\item $g_i^h:=h^{-1}g_i h$.

\item Let $N$ be a subset of $\{1,\cdots,n\}$. For each $|N|$-cycle $\al \in \fS_N$ and $\eta$ a class on $X$, 
let $\al(\eta)$ be the pullback of $\eta$ by the obvious isomorphism $X^{N}_{\al} \cong X$.

\item Two classes 
$\bigotimes_{i=1}^{\ell(\la)} g_i^{h_1}(\eta_i)$ and $\bigotimes_{i=1}^{\ell(\la)} g_i^{h_2}(\eta_i)$
on $X^n_g$ coincide for some $h_1, h_2 \in C(g)$, and a straightforward verification shows that each term
$
\bigotimes_{i=1}^{\ell(\la)} g_i^h(\eta_i)
$
repeats precisely $|\Aut(\la(\vec{\eta}))|\prod_{i=1}^{\ell(\la)} \la_i$ times. Hence, $(|\Aut(\la(\vec{\eta}))|\prod_{i=1}^{\ell(\la)} \la_i)^{-1}$ is a normalization factor to ensure that no repetition occurs in \eqref{basis}.

\item If $g=k_1 \cdots k_{\ell(\la)}$ is another cycle decomposition with $\la_i$-cycles $k_i$'s, then there exists $h\in C(g)$ such that
$$\bigotimes_{i=1}^{\ell(\la)} k_i(\eta_i)=\bigotimes_{i=1}^{\ell(\la)} g_i^h(\eta_i).$$
Thus, the expression \eqref{basis} is independent of the cycle decomposition.
\end{itemize}

Let $\fB$ be a basis for $\AT(X)$. The classes $g(\vec{\eta})$'s, with $\eta_i$'s elements of $\fB$, form a basis for $\AT(X^n_g)^{C(g)}$.

Suppose that $\hat{g}$ is another permutation of cycle type $\la$, i.e. $\hat{g}=g^\al$ for some $\al \in \fS_n$. The classes $g(\vec{\eta})$ and $\hat{g}(\vec{\eta})$ are identical in $\ATorb([\Sym^n(X)])$, though they are related by ring isomorphism 
$\al^*:\AT(X^n_{g})^{C(g)}\to \AT(X^n_{\hat{g}})^{C(\hat{g})}$ induced by $\al$. In fact, $\al^*$ sends $g(\vec{\eta})$ to $\hat{g}(\vec{\eta})$ and is independent of the choice of $\al$ due to the fact that $\vartheta^*:\AT(X^n_g)^{C(g)}\to \AT(X^n_g)^{C(g)}$ is the identity, $\forall \vartheta\in C(g)$.

We use the cohomology-weighted partition
$$\la_1(\eta_1)\cdots \la_{\ell(\la)}(\eta_{\ell(\la)}) \spa \textrm{ or simply } \spa \la(\vec{\eta})$$ 
to denote the class $g(\vec{\eta})$ (and hence $\hat{g}(\vec{\eta})$).

Now the classes $\la(\vec{\eta})$'s, running over all partitions $\la$ of $n$ and all $\eta_i \in \fB$, serve as a basis for the Chen-Ruan cohomology $\ATorb([\Sym^n(X)])$. For classes $\la(\vec{\eta})\in \ATorb([\Sym^n(X)])$ and $\rho(\vec{\xi}) \in \ATorb([\Sym^m(X)])$, keep in mind that the class
$$
\la_1(\eta_1)\cdots \la_{\ell(\la)}(\eta_{\ell(\la)})\rho_1(\xi_1)\cdots \rho_{\ell(\rho)}(\xi_{\ell(\rho)})\in \ATorb([\Sym^{n+m}(X)])
$$
is denoted by
$$
\la(\vec{\eta})\rho(\vec{\xi}).
$$
We use the shorthand $(2)$ for the divisor class $1(1)^{n-2}2(1)$. Also, we define the age of $\la(\vec{\eta})$, denote by $\age(\la(\vec{\eta}))$, to be the age of the sector $[X(\la)]$, i.e. $n-\ell(\la)$.

\subsubsection{Fixed point classes}
We can work with $\la(\vec{\eta})$'s with $\eta_k$'s in the localized cohomology $\AT(X)_\fm$ to give a basis for $\ATorb([\Sym^n(X)])_\fm$.

Assume that $X$ has exactly $p$ $\bT$-fixed points $z_1, \ldots, z_p$. For partitions $\si_1, \ldots, \si_p$, we denote the class
$$\si_{11}([z_1])\cdots \si_{1 \ell(\si_1)}([z_1]) \cdots \si_{p1}([z_p])\cdots \si_{p \ell(\si_p)}([z_p])$$ 
by $$\wt{\si}:=(\si_1, \ldots, \si_p).$$ 
The classes $\wt{\si}$'s form a basis for $\ATorb([\Sym^n(X)])_\fm$. Note also that each $\wt{\si}$ corresponds to a $\bT$-fixed point, which we denote by 
$$[\wt{\si}],$$ 
in the sector indexed by the partition 
$(\si_{11}, \ldots, \si_{1 \ell(\si_1)}, \ldots, \si_{p1}, \ldots, \si_{p \ell(\si_p)})$. So we refer to $\wt{\si}$'s as $\bT$-fixed point classes.

Moreover, given $\wt{\de} \in \ATorb([\Sym^n(X)])_\fm$ and $\wt{\si} \in \ATorb([\Sym^m(X)])_\fm$ ($m \leq n$), we say that 
$$
\wt{\de}\supset \wt{\si}
$$
if $\si_k$ is a subpartition of $\de_k$, $\forall k=1, \ldots, p$; and
$$
\wt{\de}- \wt{\si}:=(\de_1-\si_1, \ldots, \de_p-\si_p)\in \ATorb([\Sym^{n-m}(X)])_\fm.
$$
(e.g. the difference $(1,1,2,2,3)-(1,2,3)$ of two partitions is the partition $(1,2)$.)

\paragraph{$\bT$-weights.}
Given any fixed-point class $\wt{\si}$, denote by $$t(\wt{\si})$$ 
the product of $\bT$-weights on the tangent space $T_{[\wt{\si}]}\bar{I}[\Sym^m(X)]$. A simple analysis shows that
$t(\wt{\si})=\prod_{k=1}^{p} t(z_k)^{\ell(\si_k)}$. Thus, for each $\wt{\de} \supset \wt{\si}$,
$$
t(\wt{\de})=t(\wt{\si})t(\wt{\de}-\wt{\si}).
$$

When $X=\Ar$, since $L_k R_k \equiv -(r+1)^2t_1^2\mod \ttt$ for $k=1,\ldots, r+1$, it is convenient to take 
\begin{equation}\label{tau}
\tau=-(r+1)^2 t_1^2.
\end{equation}
In this manner, 
\begin{equation}\label{tangent weight}
t(\wt{\de})\equiv \tau^{\ell(\wt{\de})} \mod \ttt.
\end{equation}
Here $\ell(\wt{\de}):=\sum_{k=1}^{r+1}\ell(\de_k)$.

\paragraph{Coefficients with respect to fixed point basis.}
For $\theta({\vec{\xi}}) \in \ATorb([\Sym^m(X)])_\fm$, we have
$$\theta({\vec{\xi}})=\sum_{\wt{\si}} \frac{\langle \theta({\vec{\xi}}) |  \wt{\si} \rangle}{\langle \wt{\si} |  \wt{\si} \rangle} \wt{\si},$$
where $\langle \bullet | \bullet \rangle$ are $\bT$-equivariant orbifold pairings on $\ATorb([\Sym^m(X)])_\fm$.
Now let
$$
\al_{\theta({\vec{\xi}})}(\wt{\si}):=\frac{\langle \theta({\vec{\xi}}) |  \wt{\si} \rangle}{\langle \wt{\si} |  \wt{\si} \rangle}
$$
be the components of $\theta({\vec{\xi}})$ relative to $\wt{\si}$'s.
We have two properties by direct verification:
\begin{enumerate}
\item[(1)] Suppose $\la(\vec{\eta})$, $\rho(\vec{\varepsilon}) \in \ATorb([\Sym^n(X)])_\fm$ have explicit forms 
$\prod_{i=1}^n\prod_{j=1}^{m_i}i(\eta_{ij})$ and $\prod_{i=1}^n\prod_{j=1}^{\ell_i}i(\varepsilon_{ij})$ respectively,
we have
\[
\ \langle \la(\vec{\eta}) |  \rho(\vec{\varepsilon}) \rangle=\\
\begin{cases}
0 &\text{if $m_i\ne \ell_i$ for some $i$;}\\
\prod_{i=1}^n \langle \prod_{j=1}^{m_i}i(\eta_{ij}) \cdot \prod_{j=1}^{m_i}i(\varepsilon_{ij}) \rangle &\text{if $m_i=\ell_i$ for each $i$.}
\end{cases}
\]

\item[(2)] Given $\eta_1, \ldots, \eta_n \in \AT(X)_\fm$ and $\bT$-fixed points $y_1, \ldots, y_n$ of $X$. For $m\leq n$, the coefficient $\al_{i(\eta_1)\cdots i(\eta_n)}(i([y_1])\cdots i([y_n]))$ equals
$$\sum \al_{i(\xi_1)\cdots i(\xi_{m})}(i([y_1])\cdots i([y_m]))\al_{i(\xi_{m+1})\cdots i(\xi_{n})}(i([y_{m+1}])\cdots i([y_{n}])),$$
where the sum is over all possible $i(\xi_1)\cdots i(\xi_{m})$ and $i(\xi_{m+1})\cdots i(\xi_{n})$
such that $$i(\xi_1)\cdots i(\xi_n)=i(\eta_1)\cdots i(\eta_n).$$
\end{enumerate}
We may combine (1) with (2) to get a general statement, which presents an algorithm to calculate the coefficient $\al_{\la(\vec{\eta})}(\wt{\de})$:
\begin{prop} \label{splitting identity}
Given $\la(\vec{\eta}), \wt{\de} \in \ATorb([\Sym^n(X)])_\fm$ and $\wt{\si} \in \ATorb([\Sym^m(X)])_\fm$ with $\wt{\de}\supset \wt{\si}$,
\begin{equation}\label{coefficient}
\al_{\la(\vec{\eta})}(\wt{\de})=\sum_P  \al_{\theta({\vec{\xi}})}(\wt{\si}) \al_{\mu(\vec{\g})}(\wt{\de}-\wt{\si}),
\end{equation}
where the index $P$ under the summation symbol means that the sum is taken over all possible 
$\theta({\vec{\xi}})\in \ATorb([\Sym^m(X)])_\fm$ and $\mu(\vec{\g})\in \ATorb([\Sym^{n-m}(X)])_\fm$
satisfying $\la(\vec{\eta})=\theta({\vec{\xi}})\mu(\vec{\g})$.
\end{prop}
In the proposition, $\wt{\de}$ is separated into two parts $\wt{\si}$ and $\wt{\de}-\wt{\si}$. In general, we can break it as many parts as possible. The form \eqref{coefficient} is, however, convenient for later use.

\section{Extended Gromov-Witten theory of orbifolds}\label{s:3}
To make our exposition as self-contained as possible, we review some relevant background on orbifold Gromov-Witten theory. We take the algebro-geometric approach in the sense of Abramovich, Graber and Vistoli's works \cite{AGV1, AGV2}. The reader may also want to consult the original work \cite{ChR2} of Chen and Ruan in symplectic category.

In what follows, we utilize the isomorphism
$$A_1(\Sym^n(X); \bZ)\cong A_1(X^n; \bZ)^{\fS_n} \cong  A_1(X; \bZ).$$
In other words, we may view $E_1, \ldots, E_r$ as a basis for $A_1(\Sym^n(\Ar); \bZ)$.

\subsection{The space of twisted stable maps}
For any curve class $\be \in A_1(X; \bZ)$, the moduli space 
$$\overline{M}_{0,k}([\Sym^n(X)], \be)\footnote{ \cite{AGV1} and \cite{AGV2} adopt the notation $\mathcal{K}$ instead of $\overline{M}$. Also, we just describe the $\Spec(\bC)$-points of the moduli stack in this article. This is enough because our main purpose is the calculation of Gromov-Witten invariants.}$$ 
parametrizes genus zero $k$-pointed twisted stable map (or orbifold stable map in \cite{ChR2}) 
$$f:(\cC,\cP_1,\ldots, \cP_k)\to [\Sym^n(X)]$$ 
with the following conditions:
\begin{itemize}
\item $(\cC,\cP_1,\ldots, \cP_k)$ is an twisted nodal $k$-pointed curve. The marking $\cP_i$ is an \'etale gerbe banded by $\mu_{r_i}$, where $r_i$ is the order of the stabilizer of the twisted point. Moreover, over a node, $\cC$ has a chart isomorphic to $\Spec \spa\bC[u,v]/(uv)/\mu_{s}$ where $\mu_s$ acts on $\Spec \spa \bC[u,v]$ by $\xi\cdot (u,v) =(\xi u, \xi^{-1}v)$; locally, $c:\cC \to C$ is given by $x=u^s, y=v^s.$

\item $f$ is a representable morphism and induces a genus zero $k$-pointed stable map $f_c: (C, c(\cP_1),\ldots, c(\cP_k))\to \Sym^n(X)$ of degree $\be$ by passing to coarse moduli spaces. Note that the canonical map $c:\cC \to C$ is an isomorphism away from the nodes and marked gerbes and that whenever we say that $f$ is of degree $\be$, we actually mean $f_c$ is. 

\end{itemize}

There are evaluation maps on the moduli space $\M_{0,k}([\Sym^n(X)], \be)$, which take values in the rigidified inertia stack. At the level of $\Spec(\bC)$-points, the $i$-th evaluation map 
\begin{equation}\label{ev}
\ev_i: \overline{M}_{0,k}([\Sym^n(X)], \beta) \to \I[\Sym^n(X)]
\end{equation}
is defined by $[f:(\cC,\cP_1,\ldots, \cP_k) \to [\Sym^n(X)]] \longmapsto [f|_{\cP_i}:\cP_i \to [\Sym^n(X)]].$

The moduli space $\overline{M}_{0,k}([\Sym^n(X)], \beta)$ can be decomposed into open and closed substacks:
$$
\overline{M}_{0,k}([\Sym^n(X)], \beta)=\coprod_{\si_1,...,\si_k}\overline{M}([\Sym^n(X)],\si_1,...,\si_k; \beta).
$$
Here $\overline{M}([\Sym^n(X)],\si_1,...,\si_k; \beta)=\ev_1^{-1}([\overline{X(\si_1)}])\cap \cdots \cap \ev_k^{-1}([\overline{X(\si_k)}])$,
which can be empty for monodromy reason (e.g. the component $\overline{M}([\Sym^3(X)],(2),(2),(2); \beta)$ is empty), and the union is taken over all partitions $\si_1,...,\si_k$ of $n$. Keep in mind that the substack carries a virtual class $[\overline{M}([\Sym^n(X)],\si_1,...,\si_k; \beta)]^{\vir}$ of dimension
$$
-K_{[\Sym^n(X)]} \cdot \be + n \cdot \dim(X)+k-3-\sum_{i=1}^k \age(\si_i).
$$

The twisted map $f$ that represents an element of $\overline{M}([\Sym^n(X)],\si_1,...,\si_k; \beta)$ amounts to the following commutative diagram
\begin{equation}\label{twisted}
\begin{CD} 
P_{\cC} @> f^\prime >> X^n\\
@V VV @VV \pi V\\
\cC @>f >> [\Sym^n(X)]\\
@V c VV @VV  c V\\
C @>f_c >> \Sym^n(X)\\
\end{CD}
\end{equation}
Here $\pi$ is the natural map, $P_{\cC}:=\cC \times_{[\Sym^n(X)]} X^n$ is a scheme by representability of $f$, and $f^\prime$ is $\fS_n$-equivariant. Away from the marked points and nodes, $P_{\cC}$ is a principal $\fS_n$-bundle of $C$. It is branched over the markings with ramification types $\si_1,...,\si_k$.

Additionally, there is such a diagram
\begin{equation}\label{cover}
\begin{CD}
\tC @> \tf >> X\\
@V p  VV \\
(C, c(\cP_1), \ldots, c(\cP_k)) \\
\end{CD}
\end{equation}
associated to $f$ that $p: \tC \to C$ is an admissible cover branched over $c(\cP_1), \ldots, c(\cP_k)$ with monodromy given by $\si_1,...,\si_k$, and $\tf: \tC \to X$ is a degree $\be$ morphism such that if $\Si\subset C$ is a rational curve possessing less than $3$ special points, then there is a component of $p^{-1}(\Si)$ which is not $\tf$-contracted. In fact, \eqref{cover} is induced by the diagram \eqref{twisted} by taking $f^\prime$ mod $\fS_{n-1}$ and composing with the $n$-th projection.

The diagram \eqref{cover} will be particularly helpful later in the descriptions of $\bT$-fixed loci for the space of twisted stable maps to $[\Sym^n(\Ar)]$.
The reader should look closely at the above notation. We will use \eqref{twisted} and \eqref{cover} and the symbols there mostly without further comment. 

\subsection{Gromov-Witten invariants} \label{usual GW}
For any cohomology classes $\al_i \in \ATorb([\Sym^n(\Ar)])$ ($i=1, \ldots, k$), the $k$-point equivariant Gromov-Witten invariant is defined by
\begin{equation}\label{GW inv}
\left\langle \al_1,...,\al_k \right\rangle_{\be}^{[\Sym^n(\Ar)]}:=\int_{[\overline{M}([\Sym^n(X)], \beta)]^{\vir}_{\bT}}\ev_1^*(\al_1)\cdots \ev_k^*(\al_k),
\end{equation}
where the symbol $[\spa \spa ]^\vir_{\bT}$ stands for the $\bT$-equivariant virtual class. However, it is convenient to express the integral in \eqref{GW inv} as a sum of integrals against the virtual fundamental classes of the components $\overline{M}([\Sym^n(\Ar)],\si_1,...,\si_k; \beta)$'s.

\paragraph{Remark} The moduli space over which the integral takes is not necessarily compact. But \eqref{GW inv} is well-defined if the integral is written as a sum of residue integrals over $\bT$-fixed components via the virtual localization formula \cite{GP}. 
Alternatively, the definition \eqref{GW inv} is valid when some insertions have compact supports, e.g. $\bT$-fixed point classes. So by extending scalars, we may treat \eqref{GW inv} as a $\bQ(t_1,t_2)$-combination of invariants with at least one compactly supported insertion.
In general, the invariant takes values in $\bQ(t_1,t_2)$.
\epf

In the context of orbifolds, it is in reality more natural to study the Gromov-Witten theory in twisted degrees, i.e. in curve classes of $$A_{\orb, 1}([\Sym^n(\Ar)]; \bZ)=A_0([\overline{\Ar((2))}];\bZ)\oplus A_1([\Sym^n(\Ar)]);\bZ).$$
This makes a lot of sense because the direct sum matches $A_1(\Hilb^n(\Ar);\bZ)$ (c.f. Section \ref{s:quantum cup}). And we will see later that the $k$-point function of $[\Sym^n(\Ar)]$, with an extra quantum parameter $u$ added, has the identical number of quantum parameters to the $k$-point function of $\Hilb^n(\Ar)$ (c.f. \eqref{k-point ext} and \eqref{k-point hilb}).

Let's identify $A_0([\overline{\Ar((2))}];\bZ)$ with $\bZ$. To define the $k$-point extended Gromov-Witten invariant $\left\langle  \al_1,...,\al_k  \right\rangle_{(a,\be)}^{[\Sym^n(\Ar)]}$ of twisted degree $(a,\be) \in  \bZ \oplus A_1(\Ar;\bZ)$ with $a\geq 0$, we include additional $a$ unordered markings in the twisted stable map of degree $\be$ above such that these markings go to the age one sector under the corresponding evaluation maps. To make this precise, we present a formula:
\begin{equation} \label{GW extended}
\left\langle  \al_1,...,\al_k  \right\rangle_{(a,\be)}^{[\Sym^n(\Ar)]}=\frac{1}{a!} \left\langle  \al_1,...,\al_k, (2)^a \right\rangle_{\be}^{[\Sym^n(\Ar)]}.
\end{equation}
Note that in the expression, the last $a$ insertions are all $(2)$ and that the invariant is defined to be zero in case $a<0$. For later convenience of explanation, we refer to the markings associated to $\al_1,...,\al_k$ as distinguished marked points and to the other $a$ markings as simple marked points. Also the markings corresponding to the twisted sectors are called twisted and are otherwise called untwisted. 

The expression \eqref{GW extended} is almost identical to the non-extended version except for the appearance of the factor $\frac{1}{a!}$ due to the fact that we don't order simple markings. Additionally, we say that $\left\langle  \al_1,...,\al_k  \right\rangle_{(a,\be)}^{[\Sym^n(\Ar)]}$ is of nonzero (resp. zero) degree if it is a Gromov-Witten invariant (up to a multiple) of nonzero (resp. zero) degree and that $\left\langle  \al_1,...,\al_k  \right\rangle_{(a,\be)}^{[\Sym^n(\Ar)]}$ is multipoint if $k\geq 3$.

Like ordinary Gromov-Witten theory, if $\be \ne 0$ or $k\geq 3$, we have a forgetful morphism
$$
ft_{k+1}:\overline{M}([\Sym^n(X)],\si_1,...,\si_k, 1 ; \beta)\to \overline{M}([\Sym^n(X)],\si_1,...,\si_k; \beta)
$$
defined by forgetting the last untwisted marked points. The (untwisted) divisor equation holds as well in the orbifold case. Unfortunately, we are not allowed to forget twisted markings in general. 

\subsection{Connected version}
Let 
$$\Mc_{0,k}([\Sym^n(\Ar)], \be)$$ 
be the component of $\M_{0,k}([\Sym^n(\Ar)], \be)$ parametrizing connected covers (i.e. each cover $\tC$ associated to 
$[f:\cC \to  [\Sym^n(\Ar)]] \in \Mc_{0,k}([\Sym^n(\Ar)], \be)$
is connected).

We define $k$-point connected Gromov-Witten invariant as the contribution of the component $\Mc_{0,k}([\Sym^n(\Ar)], \be)$ to the extended Gromov-Witten invariant; namely, 
\begin{equation*}
\left\langle  \al_1, \ldots, \al_k  \right\rangle_{\be}^{[\Sym^n(\Ar)], \conn}=\int_{[\Mc_{0,k}([\Sym^n(\Ar)], \be)]^{\vir}_{\bT}}\ev_1^*(\al_1) \cdots \ev_k^*(\al_k).
\end{equation*}
Note that $\Mc_{0,k}([\Sym^n(\Ar)], \be)$ is compact whenever $\be\ne 0$, in which case the corresponding connected invariant is an element of $\bQ[t_1,t_2]$. 

Similarly, the connected invariant has an extended version. We define $k$-point extended connected invariant by
$$
\left\langle  \al_1, \ldots, \al_k  \right\rangle_{(a,\be)}^{[\Sym^n(\Ar)], \conn}=\frac{1}{a!} \left\langle  \al_1, \ldots, \al_k, (2)^a  \right\rangle_{\be}^{[\Sym^n(\Ar)], \conn}.
$$
As explained in \cite{CG}, two point extended connected invariants of $[\Sym^n(\Ar)]$ match certain connected invariants of the relative Gromov-Witten theory of $\Ar \times \bP^1$. We will explain later that the usual orbifold Gromov-Witten theory, introduced in Section \ref{usual GW}, corresponds to the relative Gromov-Witten theory with possibly disconnected source curves. 

\subsection{Orbifold quantum product}
Let $\{\omega_1,\ldots, \omega_r\}$ be the dual basis of $\{E_1,\ldots, E_r\}$ with respect to the Poincar\'e pairing. For any classes $\al_1, \ldots, \al_k \in \ATorb([\Sym^n(\Ar)])$, we define the extended $k$-point function of $[\Sym^n(\Ar)]$ by
\begin{equation}\label{k-point ext}
\llangle \al_1, \ldots, \al_k \rrangle^{[\Sym^n(\Ar)]}=\sum_{a=0}^\infty \sum_{\be \in A_1(\Ar;\bZ)} \left\langle \al_1, \ldots, \al_k \right\rangle_{(a,\beta)}^{[\Sym^n(\Ar)]}u^a s_1^{\be \cdot \omega_1 }\cdots s_r^{ \be \cdot \omega_r}
\end{equation}
and denote by
$$
\langle \al_1, \ldots, \al_k \rangle^{[\Sym^n(\Ar)]}
$$
the usual $k$-point function $\llangle \al_1, \ldots, \al_k \rrangle^{[\Sym^n(\Ar)]}|_{u=0}$. 

Now let $\{\g \}$ be a basis for the Chen-Ruan cohomology $\ATorb([\Sym^n(\Ar)])$ and $\{\g^{\vee} \}$ its dual basis. Define the small orbifold quantum product on $\ATorb([\Sym^n(\Ar)])$ in this way:
$$
\al_1 \qorb \al_2=\sum_{\g}\llangle \al_1, \al_2, \g \rrangle^{[\Sym^n(\Ar)]}\g^{\vee}.
$$
Equivalently, $\al_1 \qorb \al_2$ is defined to be the unique element satisfying
$$\langle \al_1 \qorb \al_2 \spa | \spa \al \rangle =\llangle \al_1, \al_2, \al \rrangle^{[\Sym^n(\Ar)]}, \spa \forall \al.$$
The associativity of the product follows from the WDVV equation and $1:=1(1)^n$ is the multiplicative identity because of the fundamental class axiom.
By extending scalars, we work with
$$\QATorb([\Sym^n(\Ar)])$$ 
which is defined as the vector space 
$\ATorb([\Sym^n(\Ar)])\otimes_{\bQ[t_1,t_2]}\bQ(t_1,t_2)((u, s_1, \ldots, s_r))$ endowed with quantum multiplication $\qorb$.

The extended $k$-point functions were first studied by Bryan and Graber \cite{BG} in the case of $[\Sym^n(\bC^2)]$ so as to link the Gromov-Witten theory of $[\Sym^n(\bC^2)]$ to that of $\Hilb^n(\bC^2)$. When it comes to the whole group of multipoint functions, it is clear that the extended and the usual versions share the same information. However, extended $3$-point functions are a wider group than the usual $3$-point functions, and the quantum product defined above retains more information than the usual small quantum product.

\section{Divisor operators} \label{s:4}
We are going to study the operators
$$D\qorb-$$
on the (small) quantum cohomology of the orbifold $[\Sym^n(\Ar)]$ for divisor classes $D$. We refer to them as divisor operators. We let
$$ D_k= 1(1)^{n-1}1(\omega_k), \spa k=1,\ldots, r.$$
These classes, along with $(2)$, form a basis for divisors on $[\Sym^n(\Ar)]$. Thus, the divisor operators are determined by 
$$(2)\qorb-, \spa D_1\qorb-, \ldots, D_r\qorb-,$$
which are governed by $2$-point extended invariants to be calculated in this section. 

Fix a nonnegative integer $a$ throughout the rest of this section. We shorten our notation by declaring
$$
\overline{M}([\Sym^n(\Ar)], \si_1,...,\si_k;(a,\beta))=\overline{M}([\Sym^n(\Ar)],\si_1,...,\si_k,(2)^a; \beta).
$$
Also, we use $$\vec{g}=(g_1, \cdots, g_{r+1})$$ to denote an $(r+1)$-tuple, whose entries are all partitions or all nonnegative integers.
In the case of integers, define
$$
|\vec{g}|=\ell,
$$
if the entries of $\vec{g}$ add up to $\ell$. Moreover, given a partition $\si_0$ and a multi-partition $\vec{\si}$, we put
$$\hat{\si}:=(\si_0, \vec{\si})=(\si_0, \ldots, \si_{r+1}),$$ 
which we also realize as a partition of $\sum_{k=0}^{r+1} |\si_k|$.

\subsection{Fixed loci} \label{s: fixed loci}
Let's now describe the fixed loci that will play an important role in our virtual localization calculation.

Given nonnegative integers $i$, $j$, $s$ with $1\leq i\leq j\leq r$ and $s\leq a$. We consider effective curve classes 
$$
\Eij=E_i+\cdots+E_j.
$$
For each $b^L_0 \in \{0, \ldots, s \}$ and $u^L_0 \in \{ 0, \ldots, a-s \}$, put $b^R_0=s-b^L_0$ and $u^R_0=a-s-u^L_0$. We let
\begin{equation} \label{connected fixed loci}
\{ \M_{0}^{b^L_0, \si_0, u^L_0}(1) \} \spa \spa (\textrm{resp. }\{ \M_{0}^{b^L_0, \si_0, u^L_0}(2) \})
\end{equation}
be the set consisting of all $\bT$-fixed connected components of the moduli space 
$$\Mc([\Sym^{|\la_0|}(\Ar)], \la_0, \rho_0, (2)^{s}, 1^{a-s};\spa d\Eij)$$ 
such that each point $[f_0: \cC \to [\Sym^{|\la_0|}(\Ar)]]\in \M_{0}^{b^L_0, \si_0, u^L_0}(1) \spa \spa (\textrm{resp. }\M_{0}^{b^L_0, \si_0, u^L_0}(2))$ 
has the following properties:
\begin{enumerate}
\item[(i)] $f_0$ has its source curve decomposed as
$$\cC=\cC_{L 0}\cup \cD_0 \cup \cC_{R 0}.$$ 
Here $\cC_{k0}$'s are disjoint $f_0$-contracted components, $\cD_0$ is a chain of non-contracted components with $f_{0*}([\cD_0])=d\Eij$, and $\cC_{k0}\cap \cD_0=\{ \cP_k \}$ is a twisted point, $k=L, R$.
\end{enumerate}
Let $D_0, C, P_k$ be coarse moduli spaces of $\cD_0, \cC, \cP_k$ respectively ($k=L, R$) and $\tC_0$ the admissible cover associated to $\cC$.
\begin{enumerate}
\item[(ii)] $\tC_0:=\tC_{L 0}\cup \tD_0 \cup \tC_{R 0}$ is connected with admissible covers $\tD_0 \to D_0$ and $\tC_{k0} \to C_{k0}$ ($k=L, R$). Moreover,
 \begin{itemize}
 \item each irreducible component of the cover $\tD_0 \to D_0$ is totally branched over two points (either nodes or markings) and branched nowhere else. 
 \item the covering $\tC_{L 0}\to C_{L 0}$ is branched with monodromy
 $$\la_{0}, (2)^{b^L_0}, 1^{u^L_0}, \si_0 \spa \spa (\textrm{resp. } \la_{0}, \rho_0, (2)^{b^L_0}, 1^{u^L_0}, \si_0),$$
 around markings and $P_L$ while the covering $\tC_{R 0} \to C_{R 0}$ is branched with monodromy
  $$\rho_{0}, (2)^{b^R_0}, 1^{u^R_0}, \si_0 \spa \spa (\textrm{resp. } (2)^{b^R_0}, 1^{u^R_0}, \si_0),$$
  around markings and $P_R$.
 \end{itemize} 

\item[(iii)] In the cover $\tD_0$, there exists a unique chain $\varepsilon$ formed by rational curves not contracted by $\tf_0$. Additionally,
 \begin{itemize}
 \item $\varepsilon$ possesses $j-i+1$ irreducible components which are mapped to $E_i, \ldots, E_j$ with degree $d$ under the map $\tf_0$. 
 \item the contracted components attached to the two ends of $\varepsilon$ collapse to $x_i$ and $x_{j+1}$ respectively.
 \epf
 \end{itemize}

\end{enumerate}

Now we turn our attention to the fixed locus on the moduli space 
$$\M([\Sym^n(\Ar)], \La, \wp, (a,d\Eij)).$$
We fix $\vec{b}^L$ and $\vec{b}^R$, tuples of nonnegative integers, with $|\vec{b}^L|=u^L_0$ and $|\vec{b}^R|=u^R_0$. We define
$$\cF_{\la_0, \si_0, \rho_0; b^L_0, u^L_0}^{\vec{\si}}(\vec{\la}, \vec{b}^L \spa | \spa \vec{b}^R, \vec{\rho})[i,j,s]=\{\M^{b^L_0, \si_0, u^L_0}(1)\}$$ 
to be the set of $\bT$-fixed loci of $\M([\Sym^n(\Ar)], \La, \wp, (a,d\Eij))$ (so $\La=\hat{\la}$ and $\wp=\hat{\rho}$ as partitions) with the following configuration: 
Let $[f:\cC \to [\Sym^n(\Ar)]]\in \M^{b^L_0, \si_0, u^L_0}(1)$ be a point. 

\begin{enumerate}
\item [(a)] The domain curve $\cC$ of $f$ decomposes into three pieces
\begin{equation}\label{decomposition}
\cC=\cC_L\cup \cD \cup \cC_R,
\end{equation}
where $\cC_k$'s are disjoint $f$-contracted components; $\cD$ is a chain of non-contracted components, which maps to $[\Sym^n(\Ar)]$ with degree $d\Eij$; and the intersection $\cC_k\cap \cD:=\{\cQ_k\}$ is a twisted point, $k=L, R$. 
\end{enumerate}

As in \eqref{cover}, there is an associated morphism $\tf: \tC \to \Ar$.
Let $D,C,C_k, Q_k$ be coarse moduli spaces of $\cD, \cC, \cC_k, \cQ_k$ respectively ($k=L, R$).
\begin{enumerate}
\item[(b)] $\cC_L$ carries $b^L_0+u^L_0+1$ marked points and $\cC_R$ carries the other $b^R_0+u^R_0+1=a-b^L_0-u^L_0+1$ marked points.

\item[(c)] The covering $\tC \to C$ has components 
\begin{equation}\label{covering}
\tC_k:=\tC_{L k}\cup \tD_k \cup \tC_{R k},\spa \spa k=0,\ldots, r+1.
\end{equation}
For $k\ne 0$, $\tC_k$, if nonempty, is contracted to $x_k$ in $\Ar$. (Note that $\tC_k$ is possibly empty or disconnected for $k\ne 0$. Empty sets are included just for the simplicity of notation). 

\item[(d)] For $k=0, \ldots, r+1$,
 \begin{itemize}
  \item the covering $\coprod_{k=0}^{r+1} \tC_{L k} \to C_L$ (resp. $\coprod_{k=0}^{r+1} \tC_{R k} \to C_R$) is ramified with monodromy
$$
\hat{\la}, (2)^{b^L_0+u^L_0}, \hat{\si} \spa \textrm{ (resp.  }\hat{\rho}, (2)^{b^R_0+u^R_0}, \hat{\si}),
$$
 around markings and $Q_L$ (resp. $Q_R$);
 
 \item each irreducible component of the cover $\tD_k \to D$ is totally branched over two points and branched nowhere else;
 
 \item each $\tC_{L k} \to C_L$ (resp. $\tC_{R k} \to C_R$) is a covering ramified with monodromy
 $$
 \la_{k}, (2)^{b^L_k}, 1^{b^L_0+u^L_0-b^L_k}, \si_k \spa \textrm{ (resp.  }\rho_{k}, (2)^{b^R_k}, 1^{b^R_0+u^R_0-b^R_k}, \si_k),
 $$ 
 around markings and $Q_L$ (resp. $Q_R$).
 \end{itemize}

\item[(e)] The diagram of maps
\begin{equation}\label{cover-0}
\begin{CD}
\tC_0 @> \tf|_{\tC_0} >> \Ar\\
@V   VV \\
C \\
\end{CD}
\end{equation}
corresponds to $[f_0]\in \M_{0}^{b^L_0, \si_0, u^L_0}(1)$ above.
\epf
\end{enumerate}

Note that $\cF_{\la_0, \si_0, \rho_0; b^L_0, u^L_0}^{\vec{\si}}(\vec{\la}, \vec{b}^L \spa | \spa \vec{b}^R, \vec{\rho})[i,j,s]$ does not exist for certain parameters. If it does, it is indexed by $\M_0^{b^L_0, \si_0, u^L_0}(1)$'s. Each fixed locus $\M^{b^L_0, \si_0, u^L_0}(1)$ is, however, a union of $\bT$-fixed connected components in general. 

\begin{center}
\begin{picture}(0,60)(0,0)
\put(-80,25){$\cC_L$}

\put(-50,0){\line(1,0){100}}
\put(50,0){\line(2,1){90}}
\put(-50,0){\line(-2,1){90}}
\put(-50,0){\line(-1,0){10}}
\put(50,0){\line(1,0){10}}
\put(-50,0){\line(2,-1){10}}
\put(50,0){\line(-2,-1){10}}
\put(-127,35){$\bullet$}
\put(-135,30){$\hat{\la}$}
\put(-107,25){$\bullet$}

\put(-120,20){$(2)$}

\put(-104,17){$\cdot$}
\put(-98,14){$\cdot$}
\put(-92,11){$\cdot$}

\put(-90,5){$(2)$}
\put(-77,10){$\bullet$}

\put(-53,-2){$\bullet$}
\put(-55,-10){$\hat{\si}$}

\put(48,-2){$\bullet$}
\put(50,-10){$\hat{\si}$}

\put(0,5){$\cD$}

\put(122,35){$\bullet$}
\put(130,30){$\hat{\rho}$}
\put(102,25){$\bullet$}
\put(105,20){$(2)$}

\put(87,11){$\cdot$}
\put(93,14){$\cdot$}
\put(99,17){$\cdot$}

\put(72,10){$\bullet$}
\put(75,5){$(2)$}

\put(70,25){$\cC_R$}

\end{picture}
\end{center}
\begin{quote}\small{
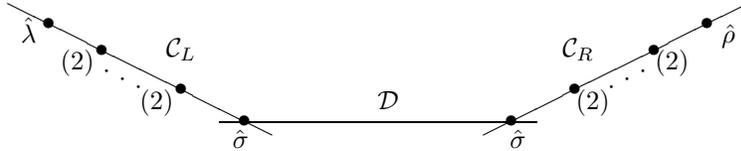
Figure 3: This is the configuration of a typical domain curve $\cC$ for the fixed locus $\M^{b^L_0, \si_0, u^L_0}(1)$. Each straight line represents a chain of curves. All markings and $\cQ_k$'s are labelled with their monodromy and there are $b^k_0+u^k_0$ copies of $(2)$ on $\cC_k$, $k=L, R$. In case $b^k_0+u^k_0=0$, $\cC_k$ is simply a twisted point. Details on the covering $\tC$ associated to $\cC$ are included in the above properties.}
\end{quote}

\medskip
Define 
$$\cF_{\la_0, \rho_0, \si_0; b^L_0, u^L_0}^{\vec{\si}}(\vec{\la},\vec{\rho}, \vec{b}^L \spa | \spa \vec{b}^R)[i,j,s]:=\{\M^{b^L_0, \si_0, u^L_0}(2)\}$$ 
in an analogous manner. The differences occur in properties (b), (d) and (e). Precisely, (b) The curve $\cC_L$ carries $b^L_0+u^L_0+2$ marked points while the curve $\cC_R$ carries the other $b^R_0+u^R_0$ marked points; (d) The covering $\coprod_{k=0}^{r+1} \tC_{L k} \to C_L$ (resp. $\coprod_{k=0}^{r+1} \tC_{R k} \to C_R$) is ramified with monodromy
$\hat{\la}, \hat{\rho}, (2)^{b^L_0+u^L_0}, \hat{\si}$ (resp. $(2)^{b^R_0+u^R_0}, \hat{\si}$) around markings and $Q_L$ (resp. $Q_R$), and the monodromy associated to the cover $\tC_{L k} \to C_L$(resp. $\tC_{R k} \to C_R$) is now $\la_{k}, \rho_{k}, (2)^{b^L_k}, 1^{b^L_0+u^L_0-b^L_k}, \si_k$ (resp. $(2)^{b^R_k}, 1^{b^R_0+u^R_0-b^R_k}, \si_k$); (e) The diagram \eqref{cover-0} corresponds to $[f_0]\in \M_{0}^{b^L_0, \si_0, u^L_0}(2)$.

\vspace{3mm}
\begin{center}
\begin{picture}(0,60)(0,0)
\put(-90,30){$\cC_L$}

\put(-50,0){\line(1,0){100}}
\put(50,0){\line(2,1){110}}
\put(-50,0){\line(-2,1){110}}
\put(-50,0){\line(-1,0){10}}
\put(50,0){\line(1,0){10}}
\put(-50,0){\line(2,-1){10}}
\put(50,0){\line(-2,-1){10}}

\put(-147,45){$\bullet$}
\put(-155,38){$\hat{\la}$}

\put(-127,35){$\bullet$}
\put(-135,30){$\hat{\rho}$}

\put(-107,25){$\bullet$}
\put(-120,20){$(2)$}
\put(-104,17){$\cdot$}
\put(-98,14){$\cdot$}
\put(-92,11){$\cdot$}
\put(-77,10){$\bullet$}
\put(-90,5){$(2)$}

\put(-53,-2){$\bullet$}
\put(-55,-10){$\hat{\si}$}

\put(48,-2){$\bullet$}
\put(50,-10){$\hat{\si}$}

\put(0,5){$\cD$}

\put(142,45){$\bullet$}
\put(148,39){$(2)$}

\put(134,30){$\cdot$}
\put(117,21){$\cdot$}
\put(99,12){$\cdot$}

\put(72,10){$\bullet$}
\put(78,5){$(2)$}

\put(80,30){$\cC_R$}

\end{picture}
\end{center}
\begin{quote}\small{
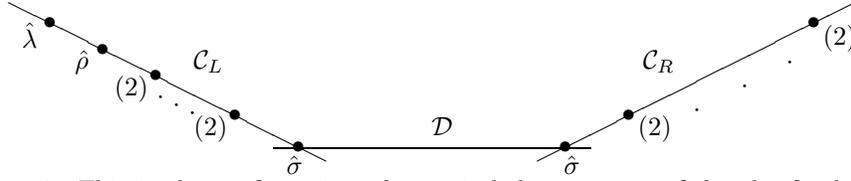
Figure 4: This is the configuration of a typical domain curve $\cC$ for the fixed locus $\M^{b^L_0, \si_0, u^L_0}(2)$. There are $b^k_0+u^k_0$ copies of $(2)$ on $\cC_k$, $k=L, R$. $\cC_L$ is always a twisted curve. $\cC_R$ is of dimension $\epsilon_2(b^R_0+u^R_0)$; in particular, it is a twisted point when $b^R_0+u^R_0 \leq 1$.}
\end{quote}

\subsection{Valuations}\label{valuation}
Given moduli space $\M([\Sym^n(\Ar)], \La, \wp, (a,d\Eij))$ as above. For each $\bT$-fixed connected component $F$, the virtual normal bundle to $F$ is denote by
$$N_F^{\vir}.$$ 
Let $[f:\cC \to [\Sym^n(\Ar)]\in F$ and $\coprod_{v}\cC_v$ the union of $1$-dimensional, contracted, connected components of $\cC$. We have a natural morphism 
$$
\phi_F: F \to F^c:=\prod_{v} \M_{0, val(v)}
$$
defined by $\phi_F([f])=([c(\cC_v)])_v$. That is, all non-contracted components, $0$-dimensional contracted components, stack structures at special points and the map $f$ are forgotten. Also, $val(v)$ denotes the number of special points on $\cC_v$. 

Let 
$$
\cF_{\la_0, \rho_0, \si_0; b^L_0, u^L_0}^{\vec{\si}}(\vec{\la}, \vec{\rho}; \vec{b}^L, \vec{b}^R)[i,j,s]
$$
be the union
$
\cF_{\la_0, \si_0, \rho_0; b^L_0, u^L_0}^{\vec{\si}}(\vec{\la}, \vec{b}^L \spa | \spa \vec{b}^R, \vec{\rho})[i,j,s]\cup \cF_{\la_0, \rho_0, \si_0; b^L_0, u^L_0}^{\vec{\si}}(\vec{\la},\vec{\rho}, \vec{b}^L \spa | \spa \vec{b}^R)[i,j,s].
$

The indices $b^L_0, \si_0, u^L_0, (k)$ ($k=L, R)$ from $\M^{b^L_0, \si_0, u^L_0}(k)$ and $\M_{0}^{b^L_0, \si_0, u^L_0}(k)$ are going to be suppressed. We simply write $\M$, $\M_{0}$.
For each $\M \in \cF_{\la_0, \rho_0, \si_0; b^L_0, u^L_0}^{\vec{\si}}(\vec{\la}, \vec{\rho}; \vec{b}^L, \vec{b}^R)[i,j,s]$, we let
$$\MT$$
be the collection of all $\bT$-fixed connected components of $\M$. 

There are other $\bT$-fixed loci on the moduli space $\M([\Sym^n(\Ar)], \La, \wp, (a,d\Eij))$. The reason why $\cF_{\la_0, \rho_0, \si_0; b^L_0, u^L_0}^{\vec{\si}}(\vec{\la}, \vec{\rho}; \vec{b}^L, \vec{b}^R)[i,j,s]$'s are singled out will be discussed later -- it will turn out that these fixed loci are enough for our study of $2$-point extended invariants as a consequence of $\ttt$-valuation below.
\begin{prop}\label{divisibility by t1+t2}
If
\begin{equation}\label{fixed loci}
\M \in \bigcup \cF_{\la_0, \rho_0, \si_0; b^L_0, u^L_0}^{\vec{\si}}(\vec{\la}, \vec{\rho}; \vec{b}^L, \vec{b}^R)[i,j,s]\spa \textrm{ and }F \in \MT,
\end{equation}
where the union ranges over all possible parameters, the inverse Euler class $\frac{1}{\eT(N_F^{\vir})}$ has valuation $1$ with respect to $\ttt$. Otherwise, it has valuation at least $2$.
\end{prop}
\bpf
Consider any $\bT$-fixed connected component $F$ of $\M([\Sym^n(\Ar)], \La, \wp ;(a,\be))$. 

Let $f:\cC \to  [\Sym^n(\Ar)]$ represent a point of $F$. As discussed earlier, there are a morphism $\tf: \tC \to \Ar$ and an ordinary stable map $f_c: C \to \Sym^n(\Ar)$ associated to $f$. Recall that $\tau=-(r+1)^2 t_1^2$. 
To establish the assertion, we need to analyze the contribution of following situations (c.f. \cite{GP}) to the inverse Euler class $\frac{1}{\eT(N_F^{\vir})}$.

\begin{enumerate}
\item \underline{Infinitesimal deformations and obstructions of $f$ with $\cC$ held fixed:}
\begin{enumerate}
\item  Any contracted component contributes zero $\ttt$-valuation. Let $\cC^\prime \subset \cC$ be a contracted component and pick any connected component $Z$ of the cover associated to $\cC^\prime$. We see that $Z$ contributes
\begin{eqnarray}\label{vertex}
\frac{\eT(H^1(Z,\tf^*T\Ar))}{\eT(H^0(Z,\tf^*T\Ar))}
\end{eqnarray}
and is collapsed by $\tf$ to $x_k$ for some $k$. So the numerator is, by Mumford's relation, congruent modulo $t_1+t_2$ to
$$
\La^{\vee}(L_k)\La^{\vee}(R_k)\equiv \tau^g,
$$
where $g=\mathrm{rank}(H^0(Z,\omega_{Z}))$ and $\La^{\vee}(t)=\sum_{i=0}^g c_{i}(H^0(Z,\omega_{Z})^{\vee}) t^{g-i}$.
The denominator of \eqref{vertex} is $\eT(T_{x_k}\Ar)$. Thus, the contribution of $Z$ is simply 
$$\tau^{g-1} \mod \ttt.$$ 
In other words, the contribution of $\cC^\prime$, being the product of the contributions of such $Z$'s, is not divisible by $t_1+t_2$.

\item The nodes joining contracted curves to non-contracted curves have zero $\ttt$-valuation because each of them gives some positive power of $\tau$ modulo $\ttt$.  

\item Non-contracted curves: Suppose $\cD$ is a non-contracted component with $\tD$ its associated (possibly disconnected) covering. Its contribution is
$$
\frac{\eT(H^1(\cD,f^*T[\Sym^n(\Ar)]))^{\mov}}{\eT(H^0(\cD,f^*T[\Sym^n(\Ar)]))^{\mov}}=\frac{\eT(H^1(\tD,\tf^*T\Ar))^{\mov}}{\eT(H^0(\tD,\tf^*T\Ar))^{\mov}}.
$$
Here $(\spa)^{\mov}$ stands for the moving part.
It is clear from (a) that each $\tf$-contracted component of $\tD$ has zero $\ttt$-valuation. However, any irreducible component $\Si$ of $\tD$ that is not $\tf$-contracted contributes
\begin{equation}\label{edge}
\frac{t_1+t_2}{\tau} \mod \ttt^2.
\end{equation}
This can be seen as follows: 

Assume that $\tf$ maps $\Si$ to $E:=\tf(\Si)$ with degree $\ell>0$. Let $S_1=\{0, \ldots 2\ell-2\}-\{\ell-1\}$ and $S_2=\{0, \ldots 2\ell\}-\{\ell\}$. 

The moving part of $\eT(H^1(\Si,\tf^*T\Ar)))$ arises from 
$$H^1(\Si, \tf^* N_{E/\Ar})=H^0(\Si, \omega_\Si \otimes \tf^* N_{E/\Ar}^\vee)^\vee.$$
The curve $E$ having self-intersection $-2$ implies $N_{E/\Ar} \cong \cO_{\bP^1}(-2)$, and so the invertible sheaf $\omega_\Si \otimes \tf^* N_{E/\Ar}^\vee$ has degree $2\ell-2$. Hence, the moving part is
$$
\ttt\prod_{k\in S_1}\frac{k(\frac{\ell-1}{\ell}(r+1)t_1)+(2\ell-2-k)(\frac{1-\ell}{\ell}(r+1)t_1)}{2\ell-2} \mod \ttt^2
$$
(which is just $\ttt$ for $\ell=1$).
We further simplify it to get
\begin{equation}\label{H1}
\ttt\tau^{\ell-1} \prod_{k=1}^{\ell-1}(\frac{\ell-k}{\ell})^2 \mod \ttt^2.
\end{equation}
On the other hand, $\eT(H^0(\Si,\tf^*T\Ar))^{\mov}$ equals $\eT(H^0(\Si, \tf^*TE))^{\mov}$, that is congruent modulo $\ttt$ to
\begin{equation}\label{H0}
\prod_{k\in S_2}\frac{k(-(r+1)t_1)+(2\ell-k)((r+1)t_1)}{2\ell} \equiv \tau^{\ell} \prod_{k=1}^{\ell-1}(\frac{\ell-k}{\ell})^2 .
\end{equation}
Dividing \eqref{H1} by \eqref{H0} gives \eqref{edge}.
\end{enumerate}

\item \underline{Infinitesimal automorphisms of $\cC$: } 
We need only investigate the non-special points of the non-contracted curves which map to fixed points. In fact, each of them gives the weight of the tangent space to the corresponding non-contracted curve and has zero $\ttt$-valuation.

\item \underline{Infinitesimal deformations of $\cC$:}

Given any node $\cP$ joining two curves $\cV_1$ and $\cV_2$. Let $P, V_1, V_2$ be coarse moduli spaces of $\cP, \cV_1, \cV_2$ respectively and $\stab(\cP)$ the stabilizer of $\cP$. In each of the following, we study the node-smoothing of $\cP$.
\begin{enumerate}
\item $\cV_1$ and $\cV_2$ are non-contracted:
We may assume that the restriction of $f_c$ to $V_k$ is a $d_k$-sheeted covering
$$f_c|_{V_k}:V_k \to \Si_k:=f_c(V_k)\cong\bP^1$$ 
for some $d_k>0$, $k=1,2.$ The node-smoothing contribution is 
\begin{equation}\label{node-smoothing 1}
|\stab(\cP)|\spa (\frac{w_1}{d_{1}}+\frac{w_2}{d_{2}})^{-1},
\end{equation}
where $w_k$ is the tangent weight of the rational curve $\Si_k$ at the fixed point $f_c(P)$.
Thus, \eqref{node-smoothing 1} is proportional to $\ttt^{-1}$ only if $d_1=d_2$ and $w_1+w_2$ is a multiple of $t_1+t_2$.

\item $\cV_1$ is non-contracted but $\cV_2$ is contracted: 
Let $w$ be the tangent weight of $V_1$ at the node $P$ and $\cL$ the tautological line bundle formed by the cotangent space $T^*_P V_2$. Denote by $\psi$ the first Chern class of $\cL$.
The node-smoothing contributes
\begin{equation}\label{node-smoothing 2}
\frac{|\stab(\cP)|}{w-\psi}.
\end{equation}
So, neither $\ttt$ nor $\ttt^{-1}$ is generated in this case.
\end{enumerate}

\end{enumerate}
Thus, only 1(c) and 3(a) may produce any power of $\ttt$.
We conclude that $F$ gives positive $\ttt$-valuation because the number of non-contracted curves is more than the number of nodes joining them. 

Suppose $F$ is a $\bT$-fixed component described in \eqref{fixed loci}, in which case we have a unique chain of non-contracted rational components for the cover associated to $\cC$. The discussion in 3(a) shows that each node in the chain gives $\ttt$-valuation $-1$. In total, the node-smoothing contributes $i-j$ in valuation. On the other hand, the chain has $j-i+1$ irreducible components. By the result of 1(c), $\frac{1}{\eT(N_F^{\vir})}$ has valuation $1$, which establishes the first assertion.

Assume that $F$ is not as in \eqref{fixed loci}. If the associated cover has at least two disjoint chains of non-contracted rational curves, a $\ttt$-valuation at least $2$ is obtained because each chain gives valuation at least $1$. Otherwise, the cover has a unique chain but property (e) (and hence (iii)) in Section \ref{s: fixed loci} is not fulfilled for each $i,j,s$. In this case, we have the same consequence by the discussion in 3(a) and the result of 1(c). This shows the second assertion.
\epf

\subsection{Counting branched covers}
Later, we will have to count certain coverings of (a chain of) rational curves. Let's now review some related notions and fix notation.

For partitions $\eta_1, \ldots, \eta_s$ of $n$, the Hurwitz number 
$$H(\eta_1, \ldots, \eta_s)$$ 
is the weighted number of possibly disconnected covers $\pi: X \to (\bP^1,p_1, \ldots, p_s)$ such that $\pi$ are branched over $p_1, \ldots, p_s$ with ramification profiles $\eta_1, \ldots, \eta_s$ and unbranched away from $p_1, \ldots, p_s$. (Each cover is counted with weight $1$ over the size of its automorphism group).

The Hurwitz number $H(\eta_1, \ldots, \eta_s)$ is essentially a combinatorial object. It can be described combinatorially by 
$$\frac{1}{n!}|\cH(\eta_1, \ldots, \eta_s)|.$$ 
Here $\cH(\eta_1, \ldots, \eta_s)$ is the set consisting of $(g_1,\ldots, g_s) \in \prod_{i=1}^s \fS_n$ satisfying (i) for each $i=1,\ldots, s$, $g_i$ has cycle type $\eta_i$; (ii) $g_1\cdots g_s=1$.

Let's introduce some other Hurwitz-type numbers.
Let 
$$\cH_{\si}(\eta_1, \ldots, \eta_s \spa|\spa \tau_1, \ldots, \tau_t)$$ 
be the subset of $\cH(\eta_1, \ldots, \eta_s, \tau_1, \ldots, \tau_t)$ such that each element $(g_1,\ldots, g_s, h_1,\ldots, h_t)$ has an additional property that $g_1\cdots g_s$ has cycle type $\si$ (and so $h_1 \cdots h_t$ has the same cycle type as well).
Put 
$$H_{\si}(\eta_1, \ldots, \eta_s \spa|\spa \tau_1, \ldots, \tau_t):=\frac{|\cH_{\si}(\eta_1, \ldots, \eta_s \spa|\spa \tau_1, \ldots, \tau_t)|}{n!}$$
(in case $\si$ is a vacuous partition, we set $H_{\si}(\eta_1, \ldots, \eta_s \spa|\spa \tau_1, \ldots, \tau_t)=1$).

We readily find the following relations:
\begin{lemma}\label{Hurwitz relation}
The number $H_{\si}(\eta_1, \ldots, \eta_s \spa|\spa \tau_1, \ldots, \tau_t)$ is exactly the product
$$|C(\si)|\spa H(\eta_1, \ldots, \eta_s, \si)\spa H(\si, \tau_1, \ldots, \tau_t).$$ 
Moreover, we have
$$
H(\eta_1,\ldots, \eta_s,\tau_1,\ldots, \tau_t)=\sum_{|\si|=n} H_{\si}(\eta_1, \ldots, \eta_s \spa|\spa \tau_1, \ldots, \tau_t).
$$
\end{lemma}

\subsection{Localization contributions}
\subsubsection{Reduction}
From now on, fix cohomology-weighted partitions $\mu_1(\vec{\eta}_1)$ and $\mu_2(\vec{\eta}_2)$ of $n$ with $\eta_{k\ell}$'s $1$ or divisors on $\Ar$. We concentrate on the $2$-point extended invariant
\begin{equation}\label{disconnected invariant}
\left\langle \mu_1(\vec{\eta}_1), \mu_2(\vec{\eta}_2)\right\rangle^{[\Sym^n(\Ar)]}_{(a, \be)}
\end{equation}
of twisted degree $(a,\be)$, $ \be \ne 0$. We will leave out the superscript $[\Sym^n(\Ar)]$ when there is no likelihood of confusion.

Let's write $$\mu_1(\vec{\eta}_1)=\kappa_1(\vec{\eta}_{11}) \theta_1(\vec{\eta}_{12}) \spa \textrm{and } \mu_2(\vec{\eta}_2)=\kappa_2(\vec{\eta}_{21}) \theta_2(\vec{\eta}_{22}),$$ 
where all entries of $\vec{\eta}_{\ell 1}$'s are 1 and all entries of $\vec{\eta}_{\ell 2}$'s are divisors, $\ell=1, 2$. We may assume that 
$$
\ell(\kappa_1)\leq \ell(\kappa_2).
$$
Use the identity
$1=\sum_{k=1}^{r+1} \frac{1}{L_k R_k}[x_k]$,
we see readily that \eqref{disconnected invariant} is a $\bQ(t_1,t_2)$-linear combination of the invariants of the form
\begin{equation}\label{fixed point}
\left\langle \kappa_{11}([x_{m_1}])\cdots \kappa_{1 \ell(\kappa_1)}([x_{m_{\ell(\kappa_1)}}])\theta_1(\vec{\eta}_{12}), \mu_2(\vec{\eta}_2)\right\rangle_{(a, \be)}.
\end{equation}
Additionally, \eqref{fixed point} is an element of $\bQ[t_1,t_2]$ as the first insertion has compact support. Also, the sum of the degrees of the insertions is at most 1 larger than the virtual dimension. Precisely, the difference is
$$
\ell(\kappa_1)-\ell(\kappa_2)+1.
$$
Thus, the invariant \eqref{fixed point} is a linear polynomial if $\ell(\kappa_1)=\ell(\kappa_2)$; otherwise, it is a rational number. 

Assume that $\be$ is not a multiple of $\Eij$ for any $i,j$. Clearly, the fixed loci \eqref{fixed loci} make no contribution. By Proposition \ref{divisibility by t1+t2}, the invariant \eqref{fixed point} is zero by divisibility of $\ttt^2$ (each of the two insertions is a linear combination of fixed-point classes with coefficients being $0$ or having nonnegative $\ttt$-valuation). It follows that \eqref{disconnected invariant} is zero as well.
So we can now set our mind on the invariant
\begin{equation}\label{disconnected invariant in dEij}
\left\langle \mu_1(\vec{\eta}_1), \mu_2(\vec{\eta}_2)\right\rangle_{(a, d\Eij)}, \spa \spa d,i,j>0.
\end{equation}

We fix positive integers $i, j, d$ with $i\leq j$ from here on. Let $\be=d\Eij$. By virtual localization, \eqref{fixed point} can be expressed as a sum of residue integrals over $\bT$-fixed loci. By Proposition \ref{divisibility by t1+t2}, the invariant \eqref{fixed point} is $\al \ttt$ for some rational number $\al$, and it suffices to evaluate \eqref{fixed point} over all $\bT$-fixed loci lying in the union
$\coprod^{i,j}\cF_{\la_0, \rho_0, \si_0; b^L_0, u^L_0}^{\vec{\si}}(\vec{\la}, \vec{\rho}; \vec{b}^L, \vec{b}^R)[i,j,s],$
where $\coprod^{i,j}$ means that only $i,j$ are fixed and the other parameters vary. Because of this, we can work modulo $\ttt^2$.

\subsubsection{Set-ups for localization}
Given $\M \in \cF_{\la_0, \rho_0, \si_0; b^L_0, u^L_0}^{\vec{\si}}(\vec{\la}, \vec{\rho}; \vec{b}^L, \vec{b}^R)[i,j,s]$ and $F \in \MT$, we let
$$
\iota_{F}: F \to \M([\Sym^n(\Ar)], \La, \wp, (a, d\Eij))
$$
be the natural inclusion (as partitions, $\La=\hat{\la}$, and $\wp=\hat{\rho}$).

Let $\M \in \cF_{\la_0, \si_0, \rho_0; b^L_0, u^L_0}^{\vec{\si}}(\vec{\la}, \vec{b}^L \spa | \spa \vec{b}^R, \vec{\rho})[i,j,s]$ (resp. $\M \in \cF_{\la_0, \rho_0, \si_0; b^L_0, u^L_0}^{\vec{\si}}(\vec{\la},\vec{\rho}, \vec{b}^L \spa | \spa \vec{b}^R)[i,j,s]$). As mentioned earlier, there are natural morphisms
$$\phi_F: F \to \M_{0,b^L_0+u^L_0+2}\times \M_{0,b^R_0+u^R_0+2} \textrm{ (resp. }\M_{0,b^L_0+u^L_0+3}\times \M_{0,b^R_0+u^R_0+1})$$ 
for $F\in \MT$ and
$$\phi_{\M_0}: \M_0 \to \M_{0,b^L_0+u^L_0+2}\times \M_{0,b^R_0+u^R_0+2} \textrm{ (resp. }\M_{0,b^L_0+u^L_0+3}\times \M_{0,b^R_0+u^R_0+1}).$$ 
Obviously, $F^c=\M_0^c$.
We intend to calculate our Gromov-Witten invariants by localization, which will be reduced to integrals over $F^c$'s. So it is necessary to understand the degree $\deg(\phi_{F})$ of the morphism $\phi_F$. 

For $F\in \MT$ with $\M \in \cF_{\la_0, \si_0, \rho_0; b^L_0, u^L_0}^{\vec{\si}}(\vec{\la}, \vec{b}^L \spa | \spa \vec{b}^R, \vec{\rho})[i,j,s]$, we consider a typical element $[f:\cC_L \cup \cD \cup \cC_R \to [\Sym^n(\Ar)]]\in F$ (in the notation of Section \ref{s: fixed loci}). The degree of $\phi_F$ is the product  $m_1 \cdot m_2$. Here 
\begin{itemize}
\item $m_1=c_0 (o(\hat{\si})^{-1}\prod_{k=1}^{r+1}|C(\si_k)|)^{\varepsilon(F)}$ is a factor arising from the nodes, which are glued over the rigidified inertia stack. Here $c_0$ is an overall factor coming from nodes of the cover $\tC_0 \to C$ (we don't have to give a careful description here as $c_0$ will be cancelled by an identical term in $\deg(\phi_{\M_0})$), and the terms $\epsilon_1(b^L_0+u^L_0)$ and $\epsilon_1(b^R_0+u^R_0)$ record the dimensions of $C_L$ and $C_R$ respectively.

\item $m_2$ is given by
$$d^{j-i+1} m_0  \prod_{k=1}^{r+1} H(\la_k, (2)^{b^L_k}, 1^{b^L_0+u^L_0-b^L_k}, \si_k) \spa H(\si_k, \si_k)^{j-i+1} \spa H(\si_k, (2)^{b^R_k}, 1^{b^R_0+u^R_0-b^R_k}, \rho_k),$$
where $d^{j-i+1}$ is an automorphism factor that takes care of the restriction $f|_{\cD}$ forgotten by $\phi_F$, $m_0$ is the contribution of $\tC_0$, and the other terms account for the overall contribution of $\coprod_{k=1}^{r+1}\tC_k$. 
\end{itemize}
Also, the degree of $\phi_{\M_0}$ can be calculated in a similar fashion. That is,
$$ 
\deg(\phi_{\M_0})= c_0 (\frac{1}{o(\si_0)})^{\varepsilon(F)} d^{j-i+1} m_0.
$$
By Lemma \ref{Hurwitz relation}, we may write $\deg(\phi_F)$ as
\begin{equation}\label{Hurwitz cover 1}
\deg(\phi_{\M_0}) (\frac{o(\si_0)}{o(\hat{\si})})^{\varepsilon(F)} \prod_{k=1}^{r+1} H_{\si_k}(\la_k, (2)^{b^L_k}, 1^{b^L_0+u^L_0-b^L_k} \spa | \spa (2)^{b^R_k}, 1^{b^R_0+u^R_0-b^R_k}, \rho_k).
\end{equation}
In the formula, $\varepsilon (F):=\epsilon_1(b^L_0+u^L_0)+\epsilon_1(b^R_0+u^R_0)+j-i$.

Similarly, for $F \in \MT$ with $\M \in \cF_{\la_0, \rho_0, \si_0; b^L_0, u^L_0}^{\vec{\si}}(\vec{\la},\vec{\rho}, \vec{b}^L \spa | \spa \vec{b}^R)[i,j,s]$, $\deg(\phi_F)$ is given by
\begin{equation}\label{Hurwitz cover 2}
\deg(\phi_{\M_0}) (\frac{o(\si_0)}{o(\hat{\si})})^{\varepsilon(F)} \prod_{k=1}^{r+1} H_{\si_k}(\la_k, \rho_k, (2)^{b^L_k}, 1^{b^L_0+u^L_0-b^L_k} \spa | \spa (2)^{b^R_k}, 1^{b^R_0+u^R_0-b^R_k}).
\end{equation}
Now $\varepsilon(F)$ is set to be $1+\epsilon_2(b^R_0+u^R_0)+j-i$.

\paragraph{Remark}
The term $(\frac{o(\si_0)}{o(\hat{\si})})^{\varepsilon(F)}$ will cancel with a similar term in $\frac{1}{\eT(N^{\vir}_{F})}$ (see Lemma \ref{euler class} below).
Moreover, forgetting the indices involving the partition $1$ does not change the value of the Hurwitz-type numbers. We did not do this in the above formulas so as to keep track of the ramification profiles corresponding to the simple marked points.

\subsubsection{Virtual normal bundles}
Let us determine $\frac{1}{\eT(N^{\vir}_{F})}$ modulo $\ttt^2$ for each connected component $F$ described in \eqref{fixed loci}.
The following outcome should be within our expectation. Recall again that $\tau=-(r+1)^2 t_1^2$.
\begin{lemma} \label{euler class}
Given $\M \in \cF_{\la_0, \rho_0, \si_0; b^L_0, u^L_0}^{\vec{\si}}(\vec{\la}, \vec{\rho}; \vec{b}^L, \vec{b}^R)[i,j,s]$ and any connected component $F \in \M_{\bT}$, we have the congruence equation
$$
\frac{1}{\eT(N^{\vir}_{F})}\equiv (\frac{o(\hat{\si})}{o(\si_0)})^{\varepsilon(F)}  \frac{\tau^{\frac{1}{2}(a-s-\ell(\vec{\la})-\ell(\vec{\rho}))}}{\eT(N^{\vir}_{\M_{0}})}  \mod \ttt^2.
$$
Here $\varepsilon(F)$'s are as in \eqref{Hurwitz cover 1}, \eqref{Hurwitz cover 2} respectively.
\end{lemma}
\bpf
We just investigate the case where $\M \in \cF_{\la_0, \si_0, \rho_0; b^L_0, u^L_0}^{\vec{\si}}(\vec{\la}, \vec{b}^L \spa | \spa \vec{b}^R, \vec{\rho})$ and $F\in \MT$, the other case being similar.

Let $p=\sum_{k=0}^{r+1}b^L_k$ and $q=\sum_{k=0}^{r+1}b^R_k$, and so $p+q=a$. Assume that $p, q > 0$.
Pick any point $[f]\in F$. Again, we follow the notation of Section \ref{s: fixed loci}. The contribution of the contracted component $\cC_L$ is
$$
\frac{\eT(H^1(\cC_L, f^*[\Sym^n(\Ar)]))}{\eT(H^0(\cC_L, f^*[\Sym^n(\Ar)]))}\equiv \tau^{\sum_k (g_k-1)}\mod \ttt.
$$
Here $g_k$'s are the genera of connected components of the covering associated to $\cC_L$. We find, by Riemann-Hurwitz formula, that
$\sum_k (g_k-1)=\frac{1}{2}(p-\ell(\hat{\la})-\ell(\hat{\si}))$.
Hence $\cC_L$ contributes
$$
\tau^{\frac{1}{2}(p-\ell(\hat{\la})-\ell(\hat{\si}))} \mod \ttt.
$$
Similarly, $\cC_R$ contributes
$$
\tau^{\frac{1}{2}(q-\ell(\hat{\rho})-\ell(\hat{\si}))} \mod \ttt.
$$
And the contribution of nodes joining contracted components to $\cD$ is 
$$\tau^{2\ell(\hat{\si})} \mod \ttt.$$ 
These three contributions, taken together, yield
$$
\tau^{\frac{1}{2}(a-\ell(\hat{\la})-\ell(\hat{\rho})+ 2\ell(\hat{\si}))} \mod \ttt.
$$
One can check that the same formula holds when $p=0$ or $q=0$.

As for the cover $\tC_{L0}\cup \tD_0 \cup \tC_{R0}$, by a similar argument, the combined contribution of $\tC_{L0}, \tC_{R0}$ and nodes joining $\tC_{L0}, \tC_{R0}$ to $\tD_0$ is given by
$$
\tau^{\frac{1}{2}(s-\ell(\la_0)-\ell(\rho_0)+ 2\ell(\si_0))} \mod \ttt.
$$
Further, the covers $\tD_1,\ldots, \tD_{r+1}$ (including the nodes inside) contribute
$$
\frac{1}{\tau^{\ell(\vec{\si})}} \mod \ttt.
$$

We now study the infinitesimal deformations of $\cC$. Let $k=L, R$. When $\cC_k$ is a curve, smoothing the node $\cP_k$ joining $\cC_k$ to $\cD$ contributes
$$\frac{o(\hat{\si})}{w_k-\psi_k},$$
where $w_k$ is the $\bT$-weight of the tangent space to $c(\cD)$ at the point $c(\cP_k)$ and $\psi_k$ is the class associated to $T^*_{c(\cP_k)} \cC_k$ (c.f. \eqref{node-smoothing 2}). By property (e) in Section \ref{s: fixed loci}, $\tf: \tC_0 \to \Ar$ corresponds to the point $[f_0: \cC_{L0}\cup \cD_0 \cup \cC_{R0} \to [\Sym^{|\la_0|}(\Ar)]] \in \M_0$, so
$$\frac{o(\si_0)}{w_k-\psi_k}$$
is the factor smoothing nodes joining $\cC_{k0}$ and $\cD_0$ and is $o(\si_0)/o(\hat{\si})$ times the preceding factor. Similarly, the overall contributions of node-smoothing inside $\cD$ and node smoothing inside $\cD_0$ differ by a factor $(o(\hat{\si})/o(\si_0))^{j-i}$. Hence, deformations of $\cC$ contribute $(o(\hat{\si})/o(\si_0))^{\varepsilon(F)}$ times those of $\cC_{L0}\cup \cD_0 \cup \cC_{R0}$, and the term
$$
(\frac{o(\hat{\si})}{o(\si_0)})^{\varepsilon(F)}\frac{1}{\eT(N^{\vir}_{\M_{0}})}
$$
is the combined contribution of the deformations of $\cC$ and the unique non-contracted connected component $\tC_0$ of the associated cover $\tC$.

Putting all these together, we get
\begin{eqnarray*}
\frac{1}{\eT(N^{\vir}_{F})}
&\equiv& (\frac{o(\hat{\si})}{o(\si_0)})^{\varepsilon(F)} \frac{1}{\eT(N^{\vir}_{\M_{0}})} \cdot \frac{\tau^{\frac{1}{2}(a-\ell(\hat{\la})-\ell(\hat{\rho})+ 2\ell(\hat{\si}))}}{\tau^{\frac{1}{2}(s-\ell(\la_0)-\ell(\rho_0)+ 2\ell(\si_0))}}\cdot \frac{1}{\tau^{\ell(\vec{\si})}}\\
&\equiv& (\frac{o(\hat{\si})}{o(\si_0)})^{\varepsilon(F)} \frac{\tau^{\frac{1}{2}(a-s-\ell(\vec{\la})-\ell(\vec{\rho}))}}{\eT(N^{\vir}_{\M_{0}})}  \mod \ttt^2,
\end{eqnarray*}
as desired.
\epf

\subsubsection{Vanishing and relation to connected invariants}
Let's look closely at the invariant \eqref{fixed point} with $\be=d\Eij$. We will go back to \eqref{disconnected invariant in dEij} in the end.

For any nonnegative integer $s$, let 
$$I(s)$$ 
be the contribution of the $\bT$-fixed loci $\coprod^{i,j,s}\cF_{\la_0, \rho_0, \si_0; b^L_0, u^L_0}^{\vec{\si}}(\vec{\la}, \vec{\rho}; \vec{b}^L, \vec{b}^R)[i,j,s]$ (all but $i,j,s$ vary) to the invariant \eqref{fixed point} with $\be=d\Eij$. We claim that
\begin{prop}\label{vanishing prop} 
For any $s<a$,
$$
I(s)\equiv 0 \mod \ttt^2.
$$
\end{prop}
Now fix a nonnegative integer $s<a$ as well. For simplicity, we drop the index $[i,j,s]$.

We would like to deduce Proposition \ref{vanishing prop} by replacing the first two insertions with $\bT$-fixed point classes. 
Fix $\bT$-fixed point classes $\vec{A}$, $\vec{B}$. Define
\begin{equation}\label{cI}
\cI:=\sum_{\M} \sum_{F \in \MT} \int_{F} \frac{\iota_{F}^*(\ev_1^*(\vec{A}) \ev_2^* (\vec{B}))}{\eT(N^{\vir}_{F})},
\end{equation}
where $\M$ is taken over all possible $\bT$-fixed loci in $\coprod_{\si_0, b^L_0, u^L_0, \vec{\si}, \vec{b}^L, \vec{b}^R} \cF_{\la_0, \rho_0, \si_0; b^L_0, u^L_0}^{\vec{\si}}(\vec{\la}, \vec{\rho}; \vec{b}^L, \vec{b}^R)$.
The coefficient 
$$
\frac{\langle \kappa_{11}([x_{m_1}])\cdots \kappa_{1 \ell(\kappa_1)}([x_{m_{\ell(\kappa_1)}}])\theta_1(\vec{\eta}_{12})  |  \vec{A}\rangle }{\langle \vec{A} |  \vec{A} \rangle} \cdot  \frac{\langle \mu_2(\vec{\eta}_2) |  \vec{B}\rangle }{\langle \vec{B} |  \vec{B} \rangle},
$$
is either zero or has nonnegative valuation with respect to $t_1+t_2$, so Proposition \ref{vanishing prop} follows from the following lemma.
\begin{lemma}\label{main lemma} 
$$
\cI \equiv 0 \mod \ttt^2.
$$
\end{lemma}
The idea of the proof is to relate $\cI$ to certain connected invariants. 

\subsubsection*{Proof of Lemma \ref{main lemma}}
Lemma \ref{main lemma} is clear if the condition 
\begin{equation}\label{condition}
\la_k\subset A_k \textrm{ and } \rho_k\subset B_k, \spa \spa \forall k=1, \ldots, r+1
\end{equation}
does not hold, in which case $\cI$ is identically zero. Now we assume \eqref{condition} and put
$$\bar{\la}_k=A_k-\la_k, \bar{\rho}_k=B_k-\rho_k.$$
That is, we may write $\vec{A}=((\la_1, \bar{\la}_1),\ldots, (\la_{r+1}, \bar{\la}_{r+1}))$ and $\vec{B}=((\rho_1, \bar{\rho}_1),\ldots, (\rho_{r+1}, \bar{\rho}_{r+1}))$.
Let
$$\bar{A}=(\bar{\la}_1, \ldots, \bar{\la}_{r+1}), \textrm{ and } \bar{B}=(\bar{\rho}_1, \ldots, \bar{\rho}_{r+1}),$$
be $\bT$-fixed point classes. 

First of all, it is good to have some observations on hand. 
\begin{lemma}\label{J}
For any partition $\si_0$ and $(r+1)$-tuples $\vec{b}^L, \vec{b}^R$, $\vec{\si}$,
$$
J_1(\si_0; b^L_0, u^L_0):= \sum_{\M \in \cF_{\la_0, \si_0, \rho_0; b^L_0, u^L_0}^{\vec{\si}}(\vec{\la}, \vec{b}^L \spa | \spa \vec{b}^R, \vec{\rho})} \deg(\phi_{\M_0}) \int_{\M_{0}^c} \frac{\iota_{\M_{0}}^*(\ev_1^*(\bar{A}) \ev_2^* (\bar{B}))}{\eT(N^{\vir}_{\M_{0}})}
$$
is 
$$ \sum_{\M \in \cF_{\la_0, \si_0, \rho_0; b^L_0, u^L_0}^{\vec{\theta}}(\vec{\la}, \vec{c}^L \spa | \spa \vec{c}^R, \vec{\rho})} \deg(\phi_{\M_0}) \int_{\M_{0}^c} \frac{\iota_{\M_{0}}^*(\ev_1^*(\bar{A}) \ev_2^* (\bar{B}))}{\eT(N^{\vir}_{\M_{0}})},
$$
and
$$
J_2(\si_0; b^L_0, u^L_0):= \sum_{\M \in \cF_{\la_0, \rho_0, \si_0; b^L_0, u^L_0}^{\vec{\si}}(\vec{\la},\vec{\rho}, \vec{b}^L \spa | \spa \vec{b}^R)} \deg(\phi_{\M_0}) \int_{\M_{0}^c} \frac{\iota_{\M_{0}}^*(\ev_1^*(\bar{A}) \ev_2^* (\bar{B}))}{\eT(N^{\vir}_{\M_{0}})}
$$
is 
$$ \sum_{\M \in \cF_{\la_0, \rho_0, \si_0; b^L_0, u^L_0}^{\vec{\theta}}(\vec{\la},\vec{\rho}, \vec{c}^L \spa | \spa \vec{c}^R)} \deg(\phi_{\M_0}) \int_{\M_{0}^c} \frac{\iota_{\M_{0}}^*(\ev_1^*(\bar{A}) \ev_2^* (\bar{B}))}{\eT(N^{\vir}_{\M_{0}})},
$$
for any $\vec{c}^L$, $\vec{c}^R$ and $\vec{\theta}$ satisfying $|\theta_k|=|\si_k|$ for each $k=1, \ldots, r+1$. Here the collections of $\bT$-fixed loci under the summation symbols are all nonempty. 
\end{lemma}
\bpf 
The first identity follows as $\cF_{\la_0, \si_0, \rho_0; b^L_0, u^L_0}^{\vec{\si}}(\vec{\la}, \vec{b}^L \spa | \spa \vec{b}^R, \vec{\rho})$ and  $\cF_{\la_0, \si_0, \rho_0; b^L_0, u^L_0}^{\vec{\theta}}(\vec{\la}, \vec{c}^L \spa | \spa \vec{c}^R, \vec{\rho})$ have the same number of elements and the same configuration for the unique non-contracted connected component of the associated cover (see the description in Section \ref{s: fixed loci}). The second identity holds for similar reasons.
\epf

We apply Proposition \ref{divisibility by t1+t2} to the connected invariant 
\begin{equation}\label{connected inv}
\left\langle \bar{A}, \bar{B}, (2)^{s}, 1^{a-s} \right\rangle_{d\Eij}^{\conn}
\end{equation}
and find that it is given by
$$
\sum_{\si_0, b^L_0, u^L_0} (J_1(\si_0; b^L_0, u^L_0)+J_2(\si_0; b^L_0, u^L_0)) \mod \ttt^2.
$$
As $a-s>0$, \eqref{connected inv} is zero. We have
\begin{equation}\label{connected invariant}
\sum_{\si_0, b^L_0, u^L_0} (J_1(\si_0; b^L_0, u^L_0)+J_2(\si_0; b^L_0, u^L_0))\equiv 0 \mod \ttt^2.
\end{equation}
 
Here is an elementary but helpful combinatorial fact.
\begin{lemma}\label{combinatorial fact}
Given nonnegative integers $k$, $p$ and $p_1, \ldots, p_k$ with $p_1+\cdots+p_k=p$. For any nonnegative integer $m \leq p$,
$$
\binom{p}{p_1, \ldots, p_k}=\sum_{m_1, \ldots, m_k}\binom{m}{m_1, \ldots, m_k}\binom{p-m}{p_1-m_1, \ldots, p_k-m_k}.
$$
Note that $\binom{\ell}{\ell_1, \ldots, \ell_k}:=0$ if $\ell$ is smaller than some of $\ell_i$'s or if some entries are negative integers.
\epf
\end{lemma}

We continue the proof of Lemma \ref{main lemma}. Let
$$
\theta=\frac{1}{2}(a-s+\ell(\vec{\la})+\ell(\vec{\rho})).
$$
For any $(r+1)$-tuple $\vec{q}$ with $|\vec{q}|=a-s$, let
$$
Q(\vec{q})=\{(\vec{b}^L, \vec{b}^R)\spa | \spa b^L_k+b^R_k=q_k, \spa \forall k=1, \ldots, r+1 \}.
$$
Fix $\si_0, b^L_0, u^L_0$, we consider two cases:
\begin{enumerate}
\item[(1)] The contribution of $\cF_{\la_0, \si_0, \rho_0; b^L_0, u^L_0}^{\vec{\si}}(\vec{\la}, \vec{b}^L \spa | \spa \vec{b}^R, \vec{\rho})$'s to $\cI$ with the constraint $(\vec{b}^L, \vec{b}^R) \in Q(\vec{q})$ is
\begin{eqnarray}\label{F}
\sum_{(\vec{b}^L, \vec{b}^R) \in Q(\vec{q})} \sum_{\vec{\si}}\sum_{\M} \sum_{F\in \MT}\int_{F} \frac{\iota_{F}^*(\ev_1^*(A) \ev_2^* (B))}{\eT(N^{\vir}_{F})} \mod \ttt^2,
\end{eqnarray}
where $\M \in \cF_{\la_0, \si_0, \rho_0; b^L_0, u^L_0}^{\vec{\si}}(\vec{\la}, \vec{b}^L \spa | \spa \vec{b}^R, \vec{\rho})$ runs through all $\bT$-fixed loci. By \eqref{tangent weight}, for each $F\in \MT$,
$$
\iota_{F}^* (\ev_1^*(A)\cdot \ev_2^*(B))\equiv \tau^{\ell(\vec{\la})+\ell(\vec{\rho})} \spa \iota_{\M_{0}}^* (\ev_1^*(\bar{A})\cdot \ev_2^*(\bar{B})) \mod \ttt.
$$ 
Applying the pushforward $\phi_{F *}$ and Lemma \ref{euler class}, \eqref{F} is given by
\begin{eqnarray*}
& &\tau^{\theta} \sum_{(\vec{b}^L, \vec{b}^R) \in Q(\vec{q})} \sum_{\vec{\si}}\sum_{\M \in \cF_{\la_0, \si_0, \rho_0; b^L_0, u^L_0}^{\vec{\si}}(\vec{\la}, \vec{b}^L \spa | \spa \vec{b}^R, \vec{\rho})} \sum_{F \in \MT} \deg(\phi_F)  \\ 
& & \cdot (\frac{o(\hat{\si})}{o(\si_0)})^{\varepsilon(F)} \int_{\M_{0}^c} \frac{\iota_{\M_{0}}^*(\ev_1^*(\bar{A}) \ev_2^* (\bar{B}))}{\eT(N^{\vir}_{\M_{0}})} \mod \ttt^2.
\end{eqnarray*}
By \eqref{Hurwitz cover 1}, \eqref{F} is congruent modulo $\ttt^2$ to
\begin{eqnarray*}
&  &\tau^{\theta}  \sum_{(\vec{b}^L, \vec{b}^R) \in Q(\vec{q})} \binom{u^L_0}{b^L_1, \ldots b^L_{r+1}}\binom{u^R_0}{b^R_1, \ldots, b^R_{r+1}}\\ 
&\cdot& \sum_{\vec{\si}} \prod_{k=1}^{r+1} H_{\si_k}(\la_k, (2)^{b^L_k}, 1^{b^L_0+u^L_0-b^L_k} \spa | \spa (2)^{b^R_k}, 1^{b^R_0+u^R_0-b^R_k}, \rho_k) \\
&  & \spa \spa \cdot \sum_{\M \in \cF_{\la_0, \si_0, \rho_0; b^L_0, u^L_0}^{\vec{\si}}(\vec{\la}, \vec{b}^L \spa | \spa \vec{b}^R, \vec{\rho})} \deg(\phi_{\M_0}) \int_{\M_{0}^c} \frac{\iota_{\M_{0}}^*(\ev_1^*(\bar{A}) \ev_2^* (\bar{B}))}{\eT(N^{\vir}_{\M_{0}})},
\end{eqnarray*}
where the product $\binom{u^L_0}{b^L_1, \ldots b^L_{r+1}}\binom{u^R_0}{b^R_1, \ldots, b^R_{r+1}}$ is the number of choices to distribute simple ramification points lying above simple markings. 

By Lemma \ref{Hurwitz relation}, \ref{J} and \ref{combinatorial fact}, \eqref{F} is simplified to
\begin{eqnarray*}
\binom{a-s}{q_1, \ldots, q_{r+1}}\tau^{\theta} \prod_{k=1}^{r+1} H(\la_k, (2)^{q_k}, 1^{a-q_k}, \rho_k) J_1(\si_0; b^L_0, u^L_0)  \mod \ttt^2.
\end{eqnarray*}

\item[(2)] By a similar argument, the contribution of $\cF_{\la_0, \rho_0, \si_0; b^L_0, u^L_0}^{\vec{\si}}(\vec{\la},\vec{\rho}, \vec{b}^L \spa | \spa \vec{b}^R)$'s to $\cI$ with the constraint $(\vec{b}^L, \vec{b}^R) \in Q(\vec{q})$ is
\begin{eqnarray*}
\binom{a-s}{q_1, \ldots, q_{r+1}}\tau^{\theta} \prod_{k=1}^{r+1} H(\la_k, (2)^{q_k}, 1^{a-q_k}, \rho_k)  J_2(\si_0; b^L_0, u^L_0) \mod \ttt^2.
\end{eqnarray*}
\end{enumerate}
In total, $\cI$ is given by
\begin{equation}\label{I}
H \cdot \sum_{\si_0, b^L_0, u^L_0}  (J_1(\si_0; b^L_0, u^L_0)+J_2(\si_0; b^L_0, u^L_0)) \mod\ttt^2,
\end{equation}
where 
$
H:=\sum_{|\vec{q}|=a-s}\binom{a-s}{q_1, \ldots, q_{r+1}}\tau^{\theta}\prod_{k=1}^{r+1} H(\la_k, (2)^{q_k}, 1^{a-q_k}, \rho_k).
$
By \eqref{connected invariant},
$$
\cI \equiv 0 \mod \ttt^2.
$$
This shows Lemma \ref{main lemma} and ends the proof of Proposition \ref{vanishing prop}.
\ep

\subsection{Combinatorial descriptions of two-point extended invariants}\label{The main result}
Now we return to the invariant \eqref{disconnected invariant in dEij}. We will interpret it combinatorially in terms of orbifold Poincar\'e pairings and connected invariants.
In general, we have the following fact on $2$-point extended invariants of nonzero degree.
\begin{theorem}\label{main result}
Given partitions $\mu_1$, $\mu_2$ of $n$ and $\ell(\mu_1)$-tuple $\vec{\eta}_1$, $\ell(\mu_2)$-tuple $\vec{\eta}_2$ with entries $1$ or divisors on $\Ar$. For any curve class $\be \ne 0$, the invariant
\begin{equation}\label{2-point}
\left\langle \mu_1(\vec{\eta}_1), \mu_2(\vec{\eta}_2) \right\rangle_{(a,\be)}
\end{equation}
is given by the sum
\begin{equation}\label{2-point sum}
\sum \langle \theta({\vec{\xi}_1})|  \theta({\vec{\xi}_2})\rangle \left\langle \nu_1(\vec{\g}_1), \nu_2(\vec{\g}_2)\right\rangle_{(a,\be)}^{\conn}.
\end{equation}
Here the sum is taken over all possible cohomology-weighted partitions $\theta({\vec{\xi}_1})$, $\theta({\vec{\xi}_2})$, $\nu_1(\vec{\g}_1)$, $\nu_2(\vec{\g}_2)$ satisfying 
$\mu_1(\vec{\eta}_1)=\theta({\vec{\xi}_1})\nu_1(\vec{\g}_1)$ and $\mu_2(\vec{\eta}_2)=\theta({\vec{\xi}_2})\nu_2(\vec{\g}_2)$.
(In particular, $\nu_1, \nu_2$ are subpartitions of $\mu_1, \mu_2$ respectively and $\mu_1-\nu_1=\theta=\mu_2-\nu_2$). 
\end{theorem}

\subsubsection*{Proof of Theorem \ref{main result}}
The statement is clear if $\be$ is not a multiple of $\Eij$ for each $i,j$ because both \eqref{2-point} and \eqref{2-point sum} vanish. Now fix $i,j,d>0$ and let $\be=d\Eij$.

We learn by Proposition \ref{vanishing prop} that $I(a)$ is the only possible contribution to \eqref{fixed point}. In other words, only 
$$\cF_{\si_0, b}(\la_0, \rho_0; \vec{\si}):=\cF_1 \cup \cF_2,$$
ranging over all possible $\la_0, \rho_0, \si_0, b, \vec{\si}$, make a contribution. Here
\begin{eqnarray}
\cF_1&=&\cF_{\la_0, \si_0, \rho_0; b, 0}^{\vec{\si}}(\vec{\si}, (0,\ldots, 0) \spa | \spa (0,\ldots, 0), \vec{\si})[i,j,a],\label{1}\\
\cF_2&=&\cF_{\la_0, \rho_0, \si_0; b, 0}^{(1^n)}(\vec{\si},\vec{\si}, (0,\ldots, 0) \spa | \spa (0,\ldots, 0))[i,j,a]\label{2}.
\end{eqnarray}
(With notation of Section \ref{s: fixed loci}, the admissible cover $\tC$ corresponding to any of these fixed loci has all its simple ramification points that are branched over simple markings in the connected component $\tC_0$ and each $\tC_k$ ($k\ne 0$) is either empty or a chain of rational curves.) 

As mentioned earlier, in order to evaluate the invariant \eqref{2-point}, it is enough to perform localization calculations over $\cF_{\si_0, b}(\la_0, \rho_0; \vec{\si})$'s because \eqref{2-point} is a linear combination of invariants of the form \eqref{fixed point}. 

We have a lemma on the inverse Euler classes of virtual normal bundles.
\begin{lemma}\label{virtual normal}
Given $F \in \MT$ with $\M \in \cF_1 \cup \cF_2$, we have
$$
\frac{1}{\eT(N^{\vir}_{F})}= (\frac{o(\hat{\si})}{o(\si_0)})^{\varepsilon_k(F)}  \frac{1}{t(\wt{\si}) \spa \eT(N^{\vir}_{\M_{0}})},
$$
for $\M \in \cF_k$, $k=1,2$. Here $\varepsilon_1(F)=\epsilon_1(b)+\epsilon_1(a-b)+j-i$ and $\varepsilon_2(F)=1+\epsilon_2(a-b)+j-i$.
\end{lemma}
\bpf
All contracted connected components of the associated cover are necessarily of genus $0$. The proof of Lemma \ref{euler class} can be carried through. 
\epf

We let $$\cI(\nu_1, \nu_2) \textrm{ and }\cI(\nu_1, \nu_2; \vec{\si})$$ 
be the contributions to \eqref{2-point} of $\coprod_{\si_0, b, \vec{\si}} \cF_{\si_0, b}(\nu_1, \nu_2; \vec{\si})$ and $\coprod_{\si_0, b} \cF_{\si_0, b}(\nu_1, \nu_2; \vec{\si})$ respectively.

Now we compute $\cI(\nu_1, \nu_2; \vec{\si})$. In order for the contribution not to vanish, the partitions $\nu_1$ and $\nu_2$ must be subpartitions of $\mu_1$ and $\mu_2$ respectively. Let us assume $\nu_1 \subset \mu_1$, $\nu_2 \subset \mu_2$. The configurations \eqref{1} and \eqref{2} force $\mu_1-\nu_1=\mu_2-\mu_2$. We set $\theta=\mu_1-\nu_1$.
\begin{lemma} \label{evaluation*}
Take a fixed locus $\M \in \coprod_{b, \si_0} \cF^{\vec{\si}}_{\si_0}(\nu_1, \nu_2, \vec{\si})$. 
For $k=1,2$ and each $F \in \MT$,
\begin{eqnarray}\label{evaluation**}
\iota_{F}^* \ev_k^*(\mu_k(\vec{\eta}_k))= t(\wt{\si}) \sum_{P_k}  \al_{\theta({\vec{\xi}_k})}(\wt{\si}) \spa \iota_{\M_{0}}^* \ev_k^*(\nu_k(\vec{\g}_k)).
\end{eqnarray}
Here $P_k$ means that we take the sum over all possible $\theta({\vec{\xi}_k})$, $\nu_k(\vec{\g}_k)$ satisfying $\mu_k(\vec{\eta}_k)=\theta({\vec{\xi}_k})\nu_k(\vec{\g}_k)$.
\end{lemma}
\bpf
The left side of \eqref{evaluation**} is $\sum_{\wt{\de}\supset \wt{\si}}\al_{\mu_k({\vec{\eta}_k})}(\wt{\de}) \spa t(\wt{\de})$.
By Proposition \ref{splitting identity}, it equals
$$
\sum_{\wt{\de}\supset \wt{\si}} \sum_{P_k} \al_{\theta({\vec{\xi}_k})}(\wt{\si}) \al_{\nu_1({\vec{\g}_k})}(\wt{\de}-\wt{\si}) \spa t(\wt{\de})=t(\wt{\si}) \sum_{P_k} \al_{\theta({\vec{\xi}_k})}(\wt{\si})\sum_{\wt{\epsilon}} \al_{\nu_1({\vec{\g}_k})}(\wt{\epsilon})\spa t(\wt{\epsilon}),
$$
which gives the right side of \eqref{evaluation**}.
\epf

It follows from Lemma \ref{evaluation*} that for each $F \in \MT$, $\iota_{F}^* (\ev_1^*(\mu_1(\vec{\eta}_1))\cdot \ev_2^*(\mu_2(\vec{\eta}_2)))$ coincides with
\begin{eqnarray}\label{evaluation}
t(\wt{\si})^2 \sum_Q \al_{\theta({\vec{\xi}_1})}(\wt{\si})  \al_{\theta({\vec{\xi}_2})}(\wt{\si}) \iota_{\M_{0}}^* (\ev_1^*(\nu_1(\vec{\g}_1))\cdot \ev_2^*(\nu_2(\vec{\g}_2))) .
\end{eqnarray}
In the formula, the index $Q$ means that the sum is over all possible $\theta({\vec{\xi}_1})$, $\theta({\vec{\xi}_2})$, $\nu_1(\vec{\g}_1))$ and $\nu_2(\vec{\g}_2))$ satisfying $\mu_1(\vec{\eta}_1)=\theta({\vec{\xi}_1})\nu_1(\vec{\g}_1)$ and  $\mu_2(\vec{\eta}_2)=\theta({\vec{\xi}_2})\nu_2(\vec{\g}_2)$.
Applying \eqref{evaluation} and Lemma \ref{virtual normal}, the contribution $\cI(\nu_1, \nu_2; \vec{\si})$ is
$$
\frac{t(\wt{\si})}{a!} \sum_Q \al_{\theta({\vec{\xi}_1})}(\wt{\si})  \al_{\theta({\vec{\xi}_2})}(\wt{\si})\sum_{\si_0, b, \M_{0}} \bH(\wt{\si}) \int_{\M_{0}} \frac{\iota_{\M_{0}}^* (\ev_1^*(\nu_1(\vec{\g}_1))\cdot \ev_2^*(\nu_2(\vec{\g}_2))}{\eT(N^{\vir}_{\M_{0}})},
$$
where $\bH(\wt{\si}):=\prod_{k=1}^{r+1} H(\si_k, \si_k)$ is a product of Hurwitz numbers. 
Thus, $\cI(\nu_1, \nu_2; \vec{\si})$ is simplified to
\begin{eqnarray*}
\bH(\wt{\si}) t(\wt{\si})  \sum_Q \al_{\theta({\vec{\xi}_1})}(\wt{\si})  \al_{\theta({\vec{\xi}_2})}(\wt{\si}) \left\langle \nu_1(\vec{\g}_1), \nu_2(\vec{\g}_2)\right\rangle_{(a,d\Eij)}^{\conn} .
\end{eqnarray*}
\\
Adding up all possible $\cI(\nu_1, \nu_2; \vec{\si})$'s, we obtain
\begin{eqnarray*}
\cI(\nu_1, \nu_2)=  \sum_Q \sum_{\wt{\si}} \bH(\wt{\si}) t(\wt{\si})  \al_{\theta({\vec{\xi}_1})}(\wt{\si})  \al_{\theta({\vec{\xi}_2})}(\wt{\si}) \left\langle \nu_1(\vec{\g}_1), \nu_2(\vec{\g}_2)\right\rangle_{(a,d\Eij)}^{\conn} .
\end{eqnarray*}
Moreover,
$$
\langle \theta({\vec{\xi}_1})|  \theta({\vec{\xi}_2})\rangle=\sum_{\wt{\si}} \al_{\theta({\vec{\xi}_1})}(\wt{\si}) \al_{\theta({\vec{\xi}_2})}(\wt{\si})\langle \wt{\si}|  \wt{\si}\rangle=\sum_{\wt{\si}}\al_{\theta({\vec{\xi}_1})}(\wt{\si})  \al_{\theta({\vec{\xi}_2})}(\wt{\si}) \bH(\wt{\si}) t(\wt{\si}).
$$
This implies that
$$
\cI(\nu_1, \nu_2)=  \sum_Q \langle \theta({\vec{\xi}_1})|  \theta({\vec{\xi}_2})\rangle \left\langle \nu_1(\vec{\g}_1), \nu_2(\vec{\g}_2)\right\rangle_{(a,d\Eij)}^{\conn}.
$$
Consequently, by taking into account of all $\cI(\nu_1, \nu_2)$'s, we deduce that \eqref{2-point} equals
$$
\sum \langle \theta({\vec{\xi}_1})|  \theta({\vec{\xi}_2})\rangle \left\langle \nu_1(\vec{\g}_1), \nu_2(\vec{\g}_2)\right\rangle_{(a,d\Eij)}^{\conn},
$$
where the sum is taken over all possible choices stated in the theorem. This finishes the proof.
\epf

The $2$-point extended connected invariants can be written explicitly in terms of certain Hurwitz numbers as in the following formulas, which are also presented in \cite{CG}.
\begin{theorem}\label{conn inv}
Let $\mu, \nu$ be partitions of $k$ and $\g_p$, $\de_q$'s 1 or divisors on $\Ar$. Let $\be \in A_1(\Ar; \bZ)$ be a nonzero curve class. Given nonnegative integer $a$, we put $g=\frac{1}{2}(a-\ell(\mu)-\ell(\nu)+2)$. If $\be=d\Eij$ for some $i,j,d>0$, the invariant $\langle \mu(\vec{\g}),\nu(\vec{\de}) \rangle_{(a,\be)}^{\conn}$ is given by
\begin{eqnarray*}
& &|\Aut(\mu)||\Aut(\nu)|\prod_{m=1}^{\ell(\mu)}(\Eij \cdot \g_m) \prod_{m=1}^{\ell(\nu)}(\Eij\cdot \de_m) \\
&\times& \frac{\ttt(-1)^{g} d^{a-1}}{k^{a-2}|\Aut(\mu(\vec{\g}))||\Aut(\nu(\vec{\de}))|}\sum_{a_1+a_2=a} \frac{H({\mu, (2)^{a_1}, (k)})\spa H({\nu, (2)^{a_2}, (k)})}{a_1!a_2!}.
\end{eqnarray*}
Otherwise, it vanishes.
\end{theorem}
\bsp
It is enough to compute
\begin{equation}\label{inv}
\left\langle \mu(\vec{\g}),\nu(\vec{\de}) \right\rangle_{(a,\be)}^{\conn}
\end{equation}
for $\g_p$, $\de_q$'s from $E_1, \ldots, E_r$ or $1$. As mentioned earlier, the invariant \eqref{inv} is a polynomial in $t_1$, $t_2$.
The second claim follows easily from Proposition \ref{divisibility by t1+t2}. To show the first claim, assume that $\be=d\Eij$ for some $i,j,d>0$. If at least one of $\g_p$, $\de_q$'s is $1$, the invariant \eqref{inv}, being a rational number, is again zero by $\ttt$-divisibility.  

Suppose that all $\g_p$, $\de_q$'s are from $E_1, \ldots, E_r$. If at least one of $\g_p$, $\de_q$'s is $E_m$ for $m\ne i,j$, \eqref{inv} also vanishes because the cover associated to the domain curve is connected and all ramification points above distinguished markings must go to either $x_i$ or $x_{j+1}$. 
Hence, it remains to evaluate \eqref{inv} with $\g_p$, $\de_q$'s from $E_i$ or $E_j$, in which case the invariant is proportional to $\ttt$.

When $r>1$, we check that $E_i \sim -\frac{1}{L_i}[x_i]$ and $E_i \sim -\frac{1}{R_{j+1}}[x_{j+1}]$ (`$\sim$' means that the difference between the left side and the right side can be written in terms of classes $[x_{i+1}], \ldots, [x_j]$ and $1$ as long as we are working modulo $\ttt$); similarly, $E_j \sim -\frac{1}{L_i}[x_i]$ and $E_j \sim -\frac{1}{R_{j+1}}[x_{j+1}]$. But when $r=1$, $E_1 \sim -\frac{2}{L_1}[x_1]$ and $E_1 \sim -\frac{2}{R_2}[x_2]$.
By the vanishing claims just verified, we can replace all $\g_p$'s with $-\frac{1}{L_i}[x_i]$ (resp. $-\frac{2}{L_1}[x_1]$) and all $\de_q$'s with $-\frac{1}{R_{j+1}}[x_{j+1}]$ (resp. $-\frac{2}{R_2}[x_2]$) for $r>1$ (resp. $r=1$). The resulting invariant is not exactly the invariant \eqref{inv}. Instead, it is congruent modulo $\ttt^2$ to 
$$
\frac{|\Aut(\mu(\vec{\g}))||\Aut(\nu(\vec{\de}))|}{|\Aut(\mu)||\Aut(\nu)|}\left\langle \mu(\vec{\g}),\nu(\vec{\de}) \right\rangle_{(a,\be)}^{\conn}.
$$
We can thus execute localization calculations over those $\bT$-fixed loci defined in \eqref{connected fixed loci} ($s=a$) by imposing one more constraint on the source curve $\cC_0$: $\cC_{L0}$ carries the marking corresponding to $\mu$ and $\cC_{R0}$ carries the marking corresponding to $\nu$ because the ramification points associated to $\mu$ (resp. $\nu$) are mapped to $x_i$ (resp. $x_{j+1}$). This means that in \eqref{connected fixed loci}, $\la_0=\mu$, $\rho_0=\nu$ and $\si_0=(k)$. In this way, we obtain the above formula in terms of Hurwitz numbers after determining precisely each term in (1)-(3) of the proof of Proposition \ref{divisibility by t1+t2}. (Here working modulo $\ttt^2$ suffices. A detailed calculation can be found in Theorem 2.3 of \cite{CG}).
\epf

The Hurwitz numbers in Theorem \ref{conn inv} can be counted explicitly. In \cite{GJV}, they are called one-part double Hurwitz numbers and satisfy
\begin{equation}\label{GJV}
\sum_b \frac{|\Aut(\si)| \spa H({\si, (2)^{b}, (k)})}{b! n^{b-1}} t^{b-\ell(\si)+1}={\frac{t/2}{sinh (t/2)}} \prod_{i=1}^{\ell(\si)}\frac{sinh (\si_it/2)}{\si_i t/2},
\end{equation}
for each partition $\si$ of $k$.
Consequently, for each nonnegative integer $b$, $H({\si, (2)^{b}, (k)})$ admits a closed-form expression.
In other words, Theorem \ref{main result} and \ref{conn inv} provide an effective method of computing $2$-point extended invariants of $[\Sym^n(\Ar)]$ of nonzero degrees. With the equations in the following proposition, this also determines the divisor operators as a consequence of $3$-point extended invariants of zero degree being determined by the Gromov-Witten theory of $[\Sym^n(\bC^2)]$. 
\begin{prop}\label{divisor eq}
Given any classes $\al_1, \ldots, \al_k \in \ATorb[\Sym^n(\Ar)]$. We have
\begin{equation}\label{(2)}
\llangle \al_1, \ldots, \al_k, (2) \rrangle= \frac{d}{du} \llangle \al_1, \ldots, \al_k \rrangle
\end{equation}
and for each $\ell=1, \ldots, r$,
\begin{equation}\label{D_l}
\llangle \al_1, \ldots, \al_k, D_\ell \rrangle= \llangle \al_1, \ldots, \al_k, D_\ell \rrangle|_{s_1, \ldots, s_r=0} + s_\ell \frac{d}{d s_\ell} \llangle \al_1, \ldots, \al_k \rrangle.
\end{equation}
\end{prop}
\bpf
By definition,
$
\left\langle \al_1, \ldots, \al_k, (2) \right\rangle_{(a,\be)}=(a+1)\left\langle \al_1, \ldots, \al_k \right\rangle_{(a+1,\be)},
$
and by the untwisted divisor equation ($\be\ne 0$ or $k \geq 3$),
$
\left\langle \al_1, \ldots, \al_k, D_\ell \right\rangle_{(a,\be)}=(\omega_\ell\cdot \be)\left\langle \al_1, \ldots, \al_k \right\rangle_{(a,\be)}.
$
These relations yield \eqref{(2)} and \eqref{D_l}. (Note, however, that \eqref{D_l} is read as
$\llangle \al_1, \ldots, \al_k, D_\ell \rrangle= s_\ell \frac{d}{d s_\ell} \llangle \al_1, \ldots, \al_k \rrangle$
for $k \geq 3$).
\epf

The sine function $\sin(u)$ is a rational function of $e^{iu}$, where $i^2=-1$. It is straightforward to verify that extended $3$-point functions involving $(2)$ or $D_\ell$ are rational functions in $t_1, t_2$, $e^{iu}, s_1, \ldots, s_r$ by the above equations.
 
On the other hand, we treat \eqref{(2)} as a twisted divisor equation because it provides a means of pulling the twisted divisor $(2)$ out. If we substitute $q=-e^{iu}$, we immediately obtain a relation on differential operators:
$$
\frac{d}{du}=i q \frac{d}{d q}.
$$
With this, \eqref{(2)} seems quite close to the usual divisor equation. They are still different, though.

\section{Comparison to other theories}\label{s:5}
\subsection{Relative Gromov-Witten theory of threefolds}\label{relative GW}
Fix $k$ distinct points $p_1, \ldots, p_k$ of $\bP^1$.
Given a positive integer $n$ and partitions $\la_1, \ldots, \la_k$ of $n$, let
$$\M^{\bullet}_{g}(\Ar \times \bP^1, (\be, n); \la_1, \ldots, \la_k)$$
be the moduli space parametrizing genus $g$ relative stable maps (c.f. \cite{L1, L2}) to $\Ar \times \bP^1$ relative to $\Ar \times p_1, \ldots, \Ar \times p_k$ with the following data:
\begin{itemize}
\item the domains are nodal curves of genus $g$ and are allowed to be disconnected;

\item the relative stable maps have degree $(\be, n)\in A_1(\Ar \times \bP^1; \bZ)$ and have nonzero degree on any connected components;  

\item the maps are ramified over the divisor $\Ar \times p_i$ with ramification type $\la_i$. The ramification points are taken to be marked and ordered.

\end{itemize}

Given any cohomology weighed partition $\la_i(\vec{\eta}_i)$, $i=1, \ldots, k$, we have an evaluation map
$$\ev_{ij}: \M^{\bullet}_{g}(\Ar \times \bP^1, (\be, n); \la_1, \ldots, \la_k)\to \Ar$$
corresponding to the ramification point of type $\la_{ij}$ over the divisor $\Ar \times p_i$. The genus $g$ relative invariant in the cohomology-weighed partitions $\la_1(\vec{\eta}_1), \ldots, \la_k(\vec{\eta}_k)$ is defined by
$$
\langle \la_1(\vec{\eta}_1), \ldots, \la_k(\vec{\eta}_k) \rangle^{\Ar \times\bP^1}_g = \frac{1}{\prod_{i=1}^k |\Aut(\la_i(\vec{\eta}_i))|} \int_{[\M^{\bullet}_{g}(\Ar \times \bP^1; \la_1, \ldots, \la_k)]_{\bT}^{\vir}}\prod_{i=1}^{k}\prod_{j=1}^{l(\la_i)} \ev_{ij}^{*}(\eta_{ij}).
$$
We define the partition function by
$$
\Zp(\Ar\times\bP^1)_{\la_1(\vec{\eta}_1), \ldots, \la_k(\vec{\eta}_k)}=\sum_g 
\langle \la_1(\vec{\eta}_1), \ldots, \la_k(\vec{\eta}_k) \rangle^{\Ar \times\bP^1}_g u^{2g-2} s_1^{\be \cdot \omega_1 }\cdots s_r^{ \be \cdot \omega_r}.
$$
However, we are more interested in the following shifted generating function
\begin{equation}\label{shifted partition}
\GW (\Ar\times\bP^1)_{\la_1(\vec{\eta}_1), \ldots, \la_k(\vec{\eta}_k)}=u^{2n-\sum_{i=1}^k \age(\la_i)} \Zp(\Ar\times\bP^1)_{\la_1(\vec{\eta}_1), \ldots, \la_k(\vec{\eta}_k)}.
\end{equation}

The partition function of relative theory of $\Ar \times \bP^1$ was studied by Maulik. We refer the reader to \cite{M} for more information. However, our results recover certain relative Gromov-Witten invariants by the following equalities. 
\begin{theorem} \label{relative theory}
For $\al=1(1)^n$, $(2)$ or $D_k$, $k=1,\ldots,r$,
\begin{equation}\label{relative inv}
\llangle \la_1(\vec{\eta}_1), \al, \la_2(\vec{\eta}_2)\rrangle^{[\Sym^n(\Ar)]}=\GW (\Ar \times \bP^1)_{\la_1(\vec{\eta}_1), \al, \la_2(\vec{\eta}_2)}.
\end{equation}
\end{theorem}
\bpf
When specialized to $s_1=\cdots=s_r=0$, the equality \eqref{relative inv} has been justified in \cite{C}. In particular, \eqref{relative inv} is valid for $\al=1(1)^n$ without the constraint.

For $\al=(2)$ or $D_k$, the coefficients of $u^i s_1^{j_1}\ldots s_r^{j_r}$, where $j_1+\cdots+j_r>0$, match up on both sides of \eqref{relative inv} by a direct comparison of Proposition 4.4 in \cite{M} with our results in Section \ref{The main result}. Hence, \eqref{relative inv} follows as well in this case.
\epf

\subsection{Quantum cohomology of Hilbert schemes of points}
\subsubsection{Nakajima basis}
The Hilbert scheme $\Hilb^n(\Ar)$ of $n$ points in $\Ar$ parametrizes $0$-dimensional closed subscheme $Z$ of $\Ar$ satisfying $$\mathrm{dim}_{\bC}H^0(Z,\cO_Z)=n.$$ 
It is a smooth irreducible variety of dimension $2n$ and provides a unique crepant resolution of $\Sym^n(\Ar)$, the so-called Hilbert-Chow morphism
$$\rho^{\HC}:\Hilb^n(\Ar) \to \Sym^n(\Ar).$$

In order to study the structure of the cohomology $\AT(\Hilb^n(\Ar)$, we review the Nakajima basis.
Given any class $\g$ in $\AT(\Ar)$, there is an Heisenberg creation operator
$$\fp_{-j}(\g):\AT(\Hilb^k(\Ar))\to A_{\bT}^{*+j-1+\mathrm{deg}(\g)}(\Hilb^{k+j}(\Ar)).$$
Let $\la$ be a partition of $n$ and $\vec{\eta}:=(\eta_1,\ldots, \eta_{\ell(\la)})$ an associated $\ell(\la)$-tuple with entries in $\AT(\Ar)_\fm$. Let $|0\rangle=1 \in A_{\bT}^0(\Ar^0)$, we define
\begin{equation}\label{Nakajima}
\fa_{\la}(\vec{\eta})=\frac{1}{|\Aut(\la(\vec{\eta}))|}\prod_{i=1}^{\ell(\la)}\frac{1}{\la_i}\spa \fp_{-\la_i}(\eta_i)|0\rangle.
\end{equation}
Choose a basis $\fB$ for $\AT(\Ar)_\fm$. The classes
$\fa_{\la}(\vec{\eta})$'s, running through all partitions $\la$ of $n$ and all $\eta_i \in \fB$,
give a basis for $\AT(\Hilb^n(\Ar))_\fm$. They are called the Nakajima basis associated to $\fB$. 

\subsubsection{Quantum cup product}\label{s:quantum cup}
Let $\rho^{\HC}_*:A_1(\Hilb^n(\Ar);\bZ) \to A_1(\Sym^n(\Ar);\bZ)$ be the homomorphism induced by the Hilbert-Chow morphism $\rho^{\HC}$. 
There are isomorphisms
$$A_1(\Hilb^n(\Ar);\bZ)\cong Ker(\rho^{\HC}_*)\oplus A_1(\Sym^n(\Ar);\bZ)\cong Ker(\rho^{\HC}_*)\oplus A_1(\Ar; \bZ).$$ 
Let $\ell$ be the class dual to the divisor $-\fa_1(1)^{n-2}a_{2}(1)$ on $\Hilb^n(\Ar)$. It is an effective rational curve class generating the kernel $Ker(\rho^{\HC}_*)$.
For any classes $\al_1, \ldots, \al_k$ on $\Hilb^n(\Ar)$, we consider the $k$-point function
\begin{equation}\label{k-point hilb}
\langle \al_1, \ldots, \al_k \rangle^{\Hilb^n(\Ar)}= \sum_{d=0}^\infty \sum_{\be \in A_1(\Ar;\bZ)}\left\langle \al_1, \ldots, \al_k \right\rangle_{(d\ell,\be)}^{\Hilb^n(\Ar)}q^d s_1^{\be \cdot \omega_1}\cdots s_r^{\be \cdot \omega_r}.
\end{equation}

Now given any basis $\{\de \}$ for $\AT(\Hilb^n(\Ar))$ and $\{\de^{\vee} \}$ its dual basis. Define the small quantum cup product $\qcrp$ on $\AT(\Hilb^n(\Ar))$ by the $3$-point functions as follows:
$$
\al_1 \qcrp \al_2=\sum_{\de}\langle \al_1, \al_2, \de \rangle^{\Hilb^n(\Ar)}\de^{\vee}
$$
Like the orbifold case, we define
$$\QAT(\Hilb^n(\Ar))$$ 
as the vector space $\AT(\Hilb^n(\Ar))\otimes_{\bQ[t_1,t_2]}\bQ(t_1,t_2)((q, s_1, \ldots, s_r))$ with the multiplication $\qcrp$. 

\subsubsection{SYM/HILB correspondence}
In Section \ref{s:4}, we provide a combinatorial description of any divisor operator on the ring $\ATorb([\Sym^n(\Ar)])$. In \cite{MO1}, on the other hand, any divisor operator on $\AT(\Hilb^n(\Ar))$ is expressed in terms of the action of affine Lie algebra $\hat{\mathfrak{gl}}(r+1)$ on the basic representations. These two expressions are actually equivalent. To make this concrete, we construct a map relating certain bases of these two cohomologies.

Let $i$ be a square root of $-1$. We make the substitution $q=-e^{iu}$. Put
$$F=\bQ(i, t_1,t_2)((u, s_1, \ldots, s_r)) \textrm{ and } K=\bQ(t_1,t_2)((u, s_1, \ldots, s_r)).$$
We define a map $L$ by
$$
L(\la(\vec{\eta}))=(-i)^{\age(\la)}\fa_{\la}(\vec{\eta}).
$$
This is obviously a one-to-one correspondence. Hence, $L$ extends to a $F$-linear isomorphism 
$$L: \QATorb([\Sym^n(\Ar)])\otimes_K F \to \QAT(\Hilb^n(\Ar))\otimes_K F.$$
Roughly speaking, $L$ maps $\bT$-fixed point basis to the Nakajima basis associated to the fixed point classes. We observe that the degree of $\fa_{\la}(\vec{\eta})$ is $n-\ell(\la)+\sum_{k=1}^{\ell(\la)} \deg(\eta_k)$, the orbifold degree of $\la(\vec{\eta})$.

Denote by $\langle \bullet | \bullet \rangle$ as well the equivariant Poincar\'e pairing on $\Hilb^n(\Ar)$.
We know that if $\la\ne \rho$, $\langle \la(\vec{\eta})|  \rho(\vec{\theta})\rangle$ and $\langle \fa_{\la}(\vec{\eta})|  \fa_{\rho}(\vec{\theta})\rangle$ vanish; otherwise, the orbifold pairing on $[\Sym^n(\Ar)]$ and the pairing on $\Hilb^n(\Ar)$ are related by
$$
\langle \la(\vec{\eta}_1) |  \la(\vec{\eta}_2)\rangle=(-1)^{\age(\la)}\langle \fa_{\la}(\vec{\eta}_1) |  \fa_{\la}(\vec{\eta}_2)\rangle.
$$
In particular, $L$ preserves (orbifold) Poincar\'e pairings, i.e. $\langle \la(\vec{\eta}) |  \rho(\vec{\theta})\rangle=\langle  L(\la(\vec{\eta})) \spa | \spa  L(\rho(\vec{\theta})) \rangle$ for all partitions $\la$, $\rho$ of $n$ and cohomology classes $\eta_i$, $\theta_j$'s of $\AT(\Ar)_\fm$.

We have the following SYM/HILB correspondence.
\begin{theorem} \label{respect divisor}
The $F$-linear isomorphism $L$ respects quantum multiplication by divisors:
\begin{equation}\label{respect divisors}
L(D \qorb \al)=L(D)\qcrp L(\al)
\end{equation}
for any class $\al$ and divisor $D$.
\end{theorem}
\bpf
For $\al=(2)$ or $D_k$ and cohomology-weighted partitions $\la_1(\vec{\eta}_1)$, $\la_2(\vec{\eta}_2)$,
\begin{eqnarray*}
\llangle \la_1(\vec{\eta}_1), \al, \la_2(\vec{\eta}_2)\rrangle^{[\Sym^n(\Ar)]}&=&\GW (\Ar \times \bP^1)_{\la_1(\vec{\eta}_1), \al, \la_2(\vec{\eta}_2)}\\
&=&\langle L(\la_1(\vec{\eta}_1)), L(\al), L(\la_2(\vec{\eta}_2)) \rangle^{\Hilb^n(\Ar)}.
\end{eqnarray*}
Indeed, the first equality is Theorem \ref{relative theory} while the second equality is Proposition 6.6 in \cite{MO1}. 

As $L$ preserves Poincar\'e pairings, it follows from the above equalities that
$$
\langle L(\la_1(\vec{\eta}_1)\qorb \al) \spa | \spa  L(\la_2(\vec{\eta}_2)) \rangle = \langle L(\la_1(\vec{\eta}_1))\qcrp L(\al) \spa | \spa L(\la_2(\vec{\eta}_2)) \rangle.
$$
This implies that $L$ respects quantum multiplication by $(2)$ and $D_k$'s. The assertion \eqref{respect divisors} now follows due to the fact that $(2)$ and $D_k$'s give a basis for divisor classes.
\epf

\subsection{An example}\label{s:example}
Let's us do some explicit calculation for the divisor operator $D_1\qorb-$ on $\ATorb([\Sym^2(\Aone)])$. 

We substitute $q=-e^{iu}$ so that
\begin{equation*}
\mathrm{sin}(\g u)=\frac{1}{2i}((-q)^{\g}-\frac{1}{\spa(-q)^{\g}}).
\end{equation*}
Consider the following basis
$$\fB:=\{1(E_1)1(E_1),\spa 
2(E_1),\spa 
1(1)1(E_1),\spa
2(1),\spa
1(1)1(1)\},$$
whose elements are ordered according to their orbifold degrees.
The matrix representation of the operator $D_1\qorb-$ with respect to $\fB$ is given by 
$$
\left(\begin{matrix}
2\theta (1-\frac{1}{1+sq}-\frac{1}{1+s/q})& i\theta(\frac{1}{1+sq}-\frac{1}{1+s/q})&-1&0&0\\
-2i\theta(\frac{1}{1+sq}-\frac{1}{1+s/q})& \theta(2-\frac{1}{1+sq}-\frac{1}{1+s/q}-\frac{2}{1-s})&0&-1&0\\
2t_1t_2&0&\frac{-\theta(1+s)}{1-s}&0&-\frac{1}{2}\\
0&4t_1t_2&0&0&0\\
0&0&4t_1t_2&0&0
\end{matrix}\right),\\$$ 
\newline 
where $\theta:=t_1+t_2$ and $s:=s_1$. This is also the matrix representation of the operator
$$L(D_1)\qcrp-$$ 
with respect to the ordered basis $L(\fB)$(c.f. \cite{MO1}). 

It is straightforward to check that $D_1 \qorb -$ has distinct eigenvalues. In particular, we have a basis $\{v_1, \ldots, v_5\}$ of eigenvectors. By quantum multiplication by $D_1$ and the identity $1$, we find 
\begin{eqnarray*}
v_i \qorb v_i&=&s_i \spa v_i, \spa \textrm{ for some } s_i\ne 0;\\
v_i \qorb v_j&=&0, \spa \spa \forall i \ne j.
\end{eqnarray*}
So by replacing $v_i$ with $v_i/s_i$, we may assume that $\{v_1, \ldots, v_5\}$ is an idempotent basis; in which case,  
\begin{equation}\label{idemp}
1=\sum_{i=1}^5 v_i.
\end{equation}

Moreover, the Vandermonde matrix associated to the eigenvalues of $D_1 \qorb -$ is invertible. In other words, by \eqref{idemp} the set 
$$\{1, D_1, D_1^2, D_1^3, D_1^4 \}$$ 
is a basis for the quantum cohomology $\QATorb([\Sym^2(\Aone)])$. Similarly, $L(D_1)$ generates the quantum ring $QA_{\bT}^*(\Hilb^2(\Aone))\otimes_K F$. We conclude that 
$$L:\QATorb([\Sym^2(\Aone)])\otimes_K F \to QA_{\bT}^*(\Hilb^2(\Aone))\otimes_K F$$ 
is indeed an $F$-algebra isomorphism.
\epf 

This simple example raises the question: Do divisor classes generate the whole quantum ring? In response to this, one may wish to examine the eigenvalues of divisor operators for bigger $n$. This, however, seems a difficult task to perform directly. 

If one of the operators $(2)\qorb -$, $D_k \qorb -$'s turns out to have distinct eigenvalues, the ring structure will be determined and $L$ will be an $F$-algebra isomorphism.
The hypothesis has yet to be entirely verified and may seems a little too good to be true. It is reasonable to expect something weaker -- maybe certain combination of these operators works. 

\section{The Crepant Resolution Conjecture} \label{s:6}
\subsection{Nonderogatory Conjecture}
We name the following nonderogatory conjecture, but we claim no originality for the statement. The reader is urged to consult \cite{MO1} for a partial evidence of the conjecture. 
\begin{conj} [\cite{MO1}] \label{nonderogatory}
Let $L$ be as in Section \ref{s:5}. The commuting family of the operators
$$L((2))\qcrp-, \spa L(D_1)\qcrp-, \ldots, L(D_r)\qcrp-$$ 
on the quantum cohomology of $\Hilb^n(\Ar)$ is nonderogatory.  That is, its joint eigenspaces are one-dimensional.
\end{conj}

Let's briefly explain some consequences of the nonderogatory conjecture on our quantum cohomology rings. Set $R=\bQ(i, t_1, t_2, q, s_1, \ldots, s_r)$ and $q=-e^{iu}$. Since the quantum ring $\AT(\Hilb^n(\Ar))\otimes_{\bQ[t_1,t_2]} R$ is semisimple, it admits a basis, say $\{v_1, \ldots, v_m\}$, of idempotent eigenvectors summing to the identity $1$. Note that the basis elements are also the simultaneous eigenvectors for $L((2))\qcrp-$,  $L(D_1)\qcrp-, \ldots, L(D_r)\qcrp-$. 

Suppose that $e_{0k}, e_{1k}, \ldots, e_{rk}$ are respectively the eigenvalues of the operators $L((2))\qcrp-$,  $L(D_1)\qcrp-, \ldots, L(D_r)\qcrp-$ corresponding to the eigenvector $v_k$. The nonderogatory property ensures that we can find numbers $a_0, a_1, \ldots, a_r$ such that 
$$
\sum_{j=0}^r a_j e_{j1}, \ldots, \sum_{j=0}^r a_j e_{jm}
$$
is a sequence of distinct elements. Therefore, the Vandermonde argument given earlier shows that the element $a_0 \cdot L((2))+\sum_{j=1}^r a_j \cdot L(D_j)$ generates $\AT(\Hilb^n(\Ar))\otimes_{\bQ[t_1,t_2]} R$. This implies that $a_0 \cdot (2)+\sum_{j=1}^r a_j \cdot D_j$ generates the quantum cohomology of $[\Sym^n(\Ar)]$ over $R$ as well. We thus obtain the following ``corollary''\footnote{ Whenever we put a double quotation mark `` '', we emphasize that the statement or word inside comes with the hypothesis of the nonderogatory conjecture.}.
\begin{conj/cor}\label{c: conjecture 1}
The divisor classes $(2)$ and $D_1, \ldots, D_r$ generate the quantum cohomology ring $\QATorb([\Sym^n(\Ar)])$, and any extended three-point function is a rational function in $t_1, t_2$, $e^{iu}, s_1, \ldots, s_r$. Under the substitution $q=-e^{iu}$, the map 
$$L:\QATorb([\Sym^n(\Ar)])\otimes_K F \to QA_{\bT}^*(\Hilb^n(\Ar))\otimes_K F$$ 
gives an isomorphism of $F$-algebras. 
\epf
\end{conj/cor}

On the other hand, we can match the orbifold Gromov-Witten theory with the relative Gromov-Witten theory:
\begin{conj/cor}
$$
\llangle \la_1(\vec{\eta}_1), \la_2(\vec{\eta}_2), \la_3(\vec{\eta}_3) \rrangle= \GW(\Ar \times \bP^1)_{\la_1(\vec{\eta}_1), \la_2(\vec{\eta}_2), \la_3(\vec{\eta}_3)},
$$
for any cohomology-weighted partitions $\la_1(\vec{\eta}_1)$, $\la_2(\vec{\eta}_2)$, $\la_3(\vec{\eta}_3)$ of $n$.
\epf
\end{conj/cor}

\subsection{Multipoint functions}
Once the nonderogatory conjecture holds, all extended $3$-point functions are known by ``Corollary'' \ref{c: conjecture 1}. In this situation, we are actually able to generalize ``Corollary'' \ref{c: conjecture 1} to cover multipoint invariants. This can be done by proceeding in an analogous manner to Okounkov and Pandharipande's determination of multipoint invariants of $\Hilb^n(\bC^2)$ (c.f. \cite{OP1}).

Let $\cB$ be a basis for the Chen-Ruan cohomology $\ATorb([\Sym^n(\Ar)])$. We recall the WDVV equation from \cite{AGV2}, but we write it in terms of extended functions to better suit our needs. For the time being, we drop the superscript $[\Sym^n(\Ar)]$. 
\begin{prop}[\cite{AGV2}]
Given Chen-Ruan classes $\al_1, \al_2, \al_3, \al_4, \be_1, \ldots, \be_k$. Let $S$ be the set $\{1, \ldots, k \}$, we have
\begin{eqnarray*}
& &\sum_{S_1 \coprod S_2=S} \sum_{\g \in \cB}\llangle \al_1, \al_2, \be_{S_1}, \g \rrangle \llangle \g^\vee, \be_{S_2}, \al_3, \al_4 \rrangle\\
&=&\sum_{S_1 \coprod S_2=S} \sum_{\g \in \cB}\llangle \al_1, \al_3, \be_{S_1}, \g \rrangle \llangle \g^\vee, \be_{S_2}, \al_2, \al_4 \rrangle. 
\end{eqnarray*}
Here, for instance, $\llangle \al_1, \al_2, \be_{S_1}, \g \rrangle:=\llangle \al_1, \al_2, \be_{i_1}, \ldots, \be_{i_\ell}, \g \rrangle$ if $S_1=\{ i_1, \ldots, i_\ell\}$.
\end{prop}

\begin{conj/prop}\label{multipoint} 
All extended multipoint functions of $[\Sym^n(\Ar)]$ can be determined from extended three-point functions and are rational functions in $t_1, t_2$, $e^{iu}, s_1, \ldots, s_r$.
\end{conj/prop}
\bpf
We may see this by induction. Suppose that any extended $m$-point function with $m\leq k$ is known and is a rational function in $t_1, t_2$, $e^{iu}, s_1, \ldots, s_r$. To determine extended $(k+1)$-point function, it suffices to study
$$
N:=\llangle \al_0, \al_1, \ldots, \al_k \rrangle
$$
for $\al_0=(2)^\ell \qorb D^{m_1}_1 \qorb\cdots \qorb D^{m_r}_r$, where $\ell$, $m_1, \ldots, m_r$ are nonnegative integers.
We may assume that $\ell+m_1+\cdots+m_r \geq 2$ in light of Proposition \ref{divisor eq} and the fundamental class axiom. Let's write
$
\al_0=D \qorb \de
$
for some $D=(2)$ or $D_j$. Clearly,
$$
N=\sum_{\g \in \cB} \llangle D, \de, \g \rrangle \llangle \g^\vee, \al_1, \ldots, \al_k \rrangle.
$$
Let $S=\{1, \ldots, k-2 \}$. By the WDVV equation,
\begin{eqnarray*}
& &\sum_{\g \in \cB} \llangle D, \de, \g \rrangle \llangle \g^\vee, \al_S, \al_{k-1}, \al_k \rrangle + \sum_{\g \in \cB} \llangle D, \de, \al_S, \g \rrangle \llangle \g^\vee, \al_{k-1}, \al_k \rrangle\\
&=&\sum_{\g \in \cB} \llangle D, \al_{k-1}, \g \rrangle \llangle \g^\vee, \al_S, \de, \al_k \rrangle + \sum_{\g \in \cB} \llangle D, \al_{k-1}, \al_S, \g \rrangle \llangle \g^\vee, \de, \al_k \rrangle \\
& & + \textrm{ (terms with extended $m$-point functions, $3 \leq m\leq k$)}.
\end{eqnarray*}
This says that $N$ is determined by lower point functions and extended $(k+1)$-point functions with a $\de$-insertion.
By replacing $D \qorb \de$ with $\de$ if necessary and continuing the above procedure, we conclude that $N$ can be calculated from lower point functions and is a rational function in $t_1, t_2$, $e^{iu}, s_1, \ldots, s_r$. By induction, our claim is thus justified.
\epf

\begin{conj/cor}[The Crepant Resolution Conjecture]
Let $q=-e^{iu}$ and $k\geq 3$. For any Chen-Ruan classes $\al_1, \ldots, \al_k$ on $[\Sym^n(\Ar)]$, we have
$$\llangle \al_1, \ldots, \al_k \rrangle^{[\Sym^n(\Ar)]}=\langle L(\al_1), \ldots, L(\al_k) \rangle^{\Hilb^n(\Ar)}.$$
In particular, $\langle \al_1, \ldots, \al_k \rrangle^{[\Sym^n(\Ar)]}=\langle L(\al_1), \ldots, L(\al_k) \rangle^{\Hilb^n(\Ar)}|_{q=-1}$.
\end{conj/cor}
\bpf
We suppress the indices $[\Sym^n(\Ar)]$ and $\Hilb^n(\Ar)$. The proof of ``Proposition'' \ref{multipoint} works as well for multipoint functions on $\Hilb^n(\Ar)$. What makes things nice is that we get exactly the same set of WDVV equations on both $[\Sym^n(\Ar)]$ and $\Hilb^n(\Ar)$ sides via $L$ provided that we have the equalities: $\llangle \al_1, \al_2, \al_3, D \rrangle=\langle L(\al_1), L(\al_2), L(\al_3), L(D) \rangle$ for $D=(2)$ and $D_j$ for $j=1, \ldots, r$. But these are clear by divisor equations and ``Corollary'' \ref{c: conjecture 1}. Thus by a recursive argument, we conclude that $L$ preserves (extended) multi-point functions and the first claim follows. The second claim is now clear.
\epf

\subsection{Closing remarks}
All ``results'' discussed above are honestly true for the case $n=2$ and $r=1$ since the divisor operator $D_1\qorb-$ has distinct eigenvalues and determines the orbifold quantum product. 

Also, in the definition of the map $L$, we may choose $-i$ instead of $i$, in which setting the correct change of variables is $q=-e^{-iu}$. As a matter of fact, the transformation $q \longmapsto \frac{1}{q}$ takes $\llangle \la_1(\vec{\eta}_1), \ldots, \la_k(\vec{\eta}_k) \rrangle^{[\Sym^n(\Ar)]}$ to $(-1)^{\sum_{j=1}^k \age(\la_j)}\llangle \la_1(\vec{\eta}_1), \ldots, \la_k(\vec{\eta}_k) \rrangle^{[\Sym^n(\Ar)]}$. To illustrate this, just look at the matrix in Section \ref{s:example}. There we observe that terms involving $q$ and $\frac{1}{q}$ agree up to a sign. 

The calculation of $[\Sym^n(\Ar)]$-invariants in Section \ref{s:4} gives an indication that these invariants might be closer, geometrically and combinatorially, to the relative invariants of $\Ar \times \bP^1$ than the invariants of $\Hilb^n(\Ar)$. In reality, it is the form the relative invariants take that motivates our calculation. 
We do know that $\GW(\Ar\times\bP^1)_{\la_1(\vec{\eta}_1), \ldots, \la_k(\vec{\eta}_k)}$ can be ``reduced'' to the $3$-point case by the degeneration formula (c.f. \cite{M}). 
It is, however, unclear if the WDVV equation ``behaves'' in a similar way to the degeneration formula. At the moment, we expect that the equality
$$
\llangle \la_1(\vec{\eta}_1), \ldots, \la_k(\vec{\eta}_k) \rrangle^{[\Sym^n(\Ar)]}=\GW(\Ar\times\bP^1)_{\la_1(\vec{\eta}_1), \ldots, \la_k(\vec{\eta}_k)}
$$
should also be true.
Particularly, the usual $k$-point function $\langle \la_1(\vec{\eta}_1), \ldots, \la_k(\vec{\eta}_k) \rangle^{[\Sym^n(\Ar)]}$ should be the coefficient of $u^{\sum_{i=1}^k \age(\la_i)-2n}$ in $\Zp(\Ar\times\bP^1)_{\la_1(\vec{\eta}_1), \ldots, \la_k(\vec{\eta}_k)}$.

\phantomsection
\addcontentsline{toc}{section}{Acknowledgements} 
\paragraph{Acknowledgements.}
This paper forms part of my Caltech Ph.D. thesis. I am especially grateful to my adviser Tom Graber for many helpful discussions and comments. 

This work is also a sequel to the joint project \cite{CG} with Amin Gholampour. It is my pleasure to thankfully acknowledge him for valuable conversations. 

Many thanks go out to Chun-Kai Li who has helped boost my morale during the preparation of this work. 

\phantomsection

\end{document}